\documentclass[11pt]{amsart}
\numberwithin{equation}{section}
\newtheorem{theorem}{Theorem}

\newtheorem{proposition}{Proposition}
\newtheorem{lemma}{Lemma}
\newtheorem{corollary}{Corollary}
\numberwithin{theorem}{section} \numberwithin{lemma}{section}
 \numberwithin{definition}{section}
\numberwithin{proposition}{section}
\numberwithin{corollary}{section}
\newtheorem{remark}{Remark}[section]
\usepackage{amsmath,amssymb}
\usepackage{mathrsfs}
\usepackage{exscale}
\usepackage{relsize}
\usepackage{color}
\usepackage{extarrows}

\def\al{\aligned}
\def\eal{\endaligned}

\def\be{\begin{equation}}
\def\ee{\end{equation}}
\def\lab{\label}

\def\t{\tilde}

\def\al{\aligned}
\numberwithin{equation}{section}

\def\q{\quad}
\def\qq{\qquad}
\def\lt{\left}
\def\rt{\right}
\def\th{\theta}
\def\p{\partial}
\def\f{\frac}
\def\nn{\nonumber}

\def\t{\tilde}

\def\s{\sqrt}
\def\la{\lambda}

\def\al{\aligned}

\def\o{\omega}
\def\O{\Omega}

\def\G{\Gamma}
\def\na{\nabla}
\def\dl{\delta}
\def\ls{\lesssim}
\def\Dl{\Delta}

\def\f{\frac}

\def\lt{\left}
\def\rt{\right}

\def\r{\rho}

\def\i{\infty}
\numberwithin{equation}{section}
\newcommand{\bR}{{\mathbb R}}

\newcommand{\tu}{\tilde{u}}
\newcommand{\bea}{\begin{eqnarray}}
\newcommand{\eea}{\end{eqnarray}}
\newcommand{\pf}{\noindent {\bf Proof. \hspace{2mm}}}
\newcommand{\ef}{ \hfill $ \Box $ \vskip 3mm}
\newcommand{\bC}{{\mathbb C}}
\newcommand{\bN}{{\mathbb N}}
\newcommand{\bZ}{{\mathbb Z}}

\newcommand{\bK}{{\bf K}}

\newcommand{\bS}{{\bf S}}

\begin{document}

\tracingpages 1
\title[D-solutions]{\bf  Decay and vanishing of some axially symmetric D-solutions of the Navier-Stokes equations}
\author{Bryan Carrillo, Xinghong Pan, Qi S. Zhang}

\address{B. Carrillo: Department of
Mathematics, University of California, Riverside, CA 92521, USA;
current address:
bcarrillo@saddleback.edu
Department of Mathematics, Saddleback College, Mission Viejo, CA 92692  }
\address{X. Pan: Department of
Mathematics, Nanjing University of Aeronautics and Astronautics, Jiang Su Province 211106, PRC }
\address{Q. S. Zhang: Department of
Mathematics, University of California, Riverside, CA 92521, USA }
\date{}

\begin{abstract}
We study axially symmetric D-solutions of the 3 dimensional Navier-Stokes equations.
The first result is an a priori decay estimate of the velocity for general domains. The second is an a priori decay estimate of the vorticity in $\bR^3$,  which improves the corresponding results in the literature. In addition, we prove a similar decay of full 3d solutions except for a small set of angles.
%, which seems to be the first quantitative decay result for general D solutions.

 Next we turn to D-solutions which are periodic in the third variable and prove vanishing result under a reasonable condition. As a corollary we prove that axially symmetric D-solutions in the slab
$\bR^2 \times I$ with suitable boundary condition is $0$. Here $I$ is any finite interval. To the best of our knowledge, this seems to be the first vanishing result on a 3 dimensional D-solution without extra  integral or decay or smallness assumption on the solution.

The tools used include Brezis-Gallouet inequality, dimension reduction, scaling, Green's function bound and Liouville theorems for Navier-Stokes equations.

\end{abstract}
\maketitle
\tableofcontents
\section{Introduction}
\q In this paper, we investigate decay and vanishing properties of axially symmetric solutions to the steady Navier-Stokes equations
\be
\left\{
\begin{aligned}
&(u\cdot\nabla)u+\nabla p-\Delta u=f, \, \text{in} \, D \subset \bR^3 \\
&\nabla\cdot u=0,
\end{aligned}
\rt. \label{1.1}
\ee
with finite Dirichlet integral
\be
\int_{D}|\nabla u(x)|^2dx<+\infty;\label{1.2}
\ee here $D \subset \bR^3$ is a noncompact, connected and axially symmetric domain on which the standard Sobolev inequality:
\[
\Vert \phi \Vert_{L^6(D)} \le S_0 \Vert \nabla \phi \Vert_{L^2(D)}
\]
 holds for all $\phi$ vanishing at infinity and for a constant $S_0$; $f$ is an divergence free
 forcing term;
it is also required that $u$ vanishes at infinity. These kind of solutions were studied in the pioneer work of Leray \cite{Le:1} (p24) by variational method and are often referred as D-solutions.
If $f=0$ and $u=0$ on $\p D$, then the solution (not necessarily axially symmetric) is called a homogeneous D-solution. The following has been an  old open question:

 \noindent {\it Is a homogeneous D-solution equal to  $0$ ?}

Despite the apparent simplicity, it is not even known if a general D-solution has any definite decay rate comparing with the distance function near infinity, even when the domain is $\bR^3$.
Under some extra integral or decay assumptions for the solution $u$, vanishing results were derived by a number of authors.  For instance, Galdi \cite{Ga:1} Theorem X.9.5 proved that if $u$ is a homogeneous D-solution in the domain $D=\bR^3$ and $u \in L^{9/2}(\bR^3)$, then $u=0$.  This result was improved by a log factor in Chae and Wolf \cite{CW:1}. In the paper \cite{Ch:1}, Chae proved that homogeneous D-solutions in $\bR^3$ is $0$ if also $\Delta u \in L^{6/5}(\bR^3)$, a condition scaled the same way as $\Vert \nabla u \Vert_2$. Seregin \cite{Se:1} proved that homogeneous D-solutions in $\bR^3$ is $0$  if $u \in L^6(\bR^3) \cap BMO^{-1}$. In a recent paper \cite{KTW:1}, Kozono etc showed that if the vorticity $w=w(x)$ decays faster than $C/|x|^{5/3}$ at infinity, then
homogeneous D-solutions in $\bR^3$ is $0$. Under certain smallness assumption, vanishing result for homogeneous 3 dimensional solutions in a slab was also obtained in the book \cite{Ga:1}, Chapter XII.

We will concentrate on the axially symmetric homogeneous D-solutions in this paper. The following
 cylindrical coordinates will be used through out: $x=(x_1, x_2, z)$, $\th= \tan^{-1}{x_2/x_1}$  and $r=\sqrt{x^2_1+x^2_2}$. $e_r =(x_1/r, x_2/r, 0)$, $e_\th=(-x_2/r, x_1/r, 0)$ and $e_z=(0, 0, 1)$. For convenience, we often write $x'=(x_1, x_2)$ and $x_3=z$.

Recall that if a smooth vector field $u(x)=u^r(r,z)e_r+u^\th(r,z)e_\th+u^z(r,z)e_z$ is an axially symmetric solution of \eqref{1.1}, then $u^r, u^\th, u^z$ satisfy the following equations
\be
\lt\{
\begin{aligned}
&(u^r\p_r+u^z\p_z)u^r-\f{(u^\th)^2}{r}+\p_r p=(\p^2_r+\f{1}{r}\p_r+\p^2_z-\f{1}{r^2})u^r,\\
&(u^r\p_r+u^z\p_z)u^\th+\f{u^ru^\th}{r}=(\p^2_r+\f{1}{r}\p_r+\p^2_z-\f{1}{r^2})u^\th,\\
&(u^r\p_r+u^z\p_z)u^z+\p_z p=(\p^2_r+\f{1}{r}\p_r+\p^2_z)u^z,\\
&\p_ru^r+\f{u^r}{r}+\p_zu^z=0.
\end{aligned}
\rt.
 \label{1.4}
\ee

The vorticity  $w$ is defined as $w(x)=\nabla\times u(x)=w^r(r,z)e_r+w^\th(r,z)e_\th+w^z(r,z)e_z$, where
\be
w^r=-\p_zu^\th,\q w^\th=\p_zu^r-\p_ru^z, \q w^z=\f{1}{r}\p_r(ru^\th). \nn
\ee
The equations for $w^r, w^\th, w^z$ are
\be
\lt\{
\begin{aligned}
&(u^r\p_r+u^z\p_z)w^r-(w^r\p_r+w^z\p_z)u^r=(\p^2_r+\f{1}{r}\p_r+\p^2_z-\f{1}{r^2})w^r,\\
&(u^r\p_r+u^z\p_z)w^\th-\f{u^r}{r}w^\th-\f{1}{r}\p_z(u^\th)^2=(\p^2_r+\f{1}{r}\p_r+\p^2_z-\f{1}{r^2})w^\th,\\
&(u^r\p_r+u^z\p_z)w^z-(w^r\p_r+w^z\p_z)u^z=(\p^2_r+\f{1}{r}\p_r+\p^2_z)w^z.\label{1.5}
\end{aligned}
\rt.
\ee

 Although the solution $u$ is independent of the angle $\theta$ in the cylindrical system, the aforementioned  question is also wide open. However certain decay estimates exist for the homogeneous D-solution $u$ and vorticity $w$.
More specifically, the combined result of Chae-Jin \cite{CJ:1} and Weng \cite{We:1} state that, for $x \in \bR^3$,
\be
\lab{cjw}
|u(x)| \le C \left(\frac{\log r}{r} \right)^{1/2}, \q |w^\th(x)| \le C r^{-(19/16)^-}, \q
|w^r(x)| + |w^z(x)| \le C r^{-(67/64)^-}.
\ee Here $C$ is a positive constant and for a positive number $a$, we write $a^-$ as a number which is smaller than but close to $a$. These authors use line integral techniques from the work on 2 dimensional D-solutions by Gilbarg and Weinberger \cite{GW:1}. In a recent paper \cite{LZ:1}, among other things, Liouville property for bounded,  axially symmetric solutions of the Navier-Stokes equation were studied.
Assuming in addition that $r u^\theta$ is bounded and $u$ is periodic in the $z$ variable, then
it was shown that $u \equiv 0$.

 The first result of this paper is an a priori decay estimate of the velocity for general domains. The second result is an a priori decay estimate of the vorticity in $\bR^3$,  which improves the corresponding results \eqref{cjw} in the literature. In the process of the proof, we also give a very short proof of the decay estimate of the velocity in \eqref{cjw} under a more general condition. Next we turn to D-solutions which are periodic in the third variable and prove vanishing result under a reasonable condition. As a corollary we prove that axially symmetric D-solutions in
$\bR^2 \times I$ with suitable boundary condition is $0$. Here $I$ is any finite interval. To the best of our knowledge, this seems to be the first vanishing result on D-solutions without extra integral or decay assumption on the solution.

Now let us present the main results in detail. In most of this paper, typically the domain $D$ is either $\bR^3$ or $\bR^2 \times S^1$. In the later case the flow is periodic in the $Z$ direction with period $2 \pi$, a number chosen for convenience. Any other positive period is fine. We will always take $f=0$ through out the paper, although the decay estimates are still valid for $f$ that decays sufficiently fast.

The next three theorems are the main results of the paper, which cover a little more general class of solutions, namely the Dirichlet integral is allowed to be log divergent.

\begin{theorem}
\lab{thfullwdecay}
Let $u$ be a smooth axially symmetric solution to the problem
\be
\left\{
\begin{aligned}
&(u\cdot\nabla)u+\nabla p-\Delta u=0, \quad \text{in} \quad  \bR^3,\\
&\nabla\cdot u=0,\\
&\lim\limits_{|x|\rightarrow \i}u=0,
\end{aligned}
\rt.
\ee
such that the Dirichlet integral satisfies the condition: for a constant $C$, and all $R \ge 1$,
\be\label{dr2r1}
\int^\i_{-\i}\int_{R \le |x'| \le 2 R}\big(|\nabla u(x)|^2 dx+ |u(x)|^6\big)dx'dx_3<C<\infty.
\ee
 Then the velocity and vorticity satisfy the following a priori bound. For a constant $C_0>0$, depending only on the constant $C$ in (\ref{dr2r1}) such that
 \[
 |u(x)| \le C_0 \frac{(\ln r)^{1/2}}{\sqrt{r}};
 \]
\[
\al
&|w^\th(x)| \le C_0  \f{(\ln r)^{3/4}}{r^{5/4}}, \quad |w^r(x)|+|w^z(x)| \le C \f{(\ln r)^{11/8}}{r^{9/8}}, \q r \ge e.
\eal
\]

\end{theorem}

\begin{theorem}
\lab{thzp}
Let $u$ be a smooth axially symmetric solution to the problem
\be\label{e10.2}
\left\{
\begin{aligned}
&(u\cdot\nabla)u+\nabla p-\Delta u=0, \quad \text{in} \quad  \bR^2 \times S^1=\bR^2 \times [-\pi, \pi],\\
&\nabla\cdot u=0,\\
&u(x_1,x_2,z)=u(x_1,x_2,z+2\pi),\\
&\lim\limits_{|x|\rightarrow \i}u=0,
\end{aligned}
\rt.
\ee
such that the Dirichlet integral satisfies the condition: for a constant $C$, and all $R \ge 1$,
\be
\lab{dr2r}
\int^\pi_{-\pi}\int_{R \le |x'| \le 2 R}|\nabla u(x)|^2 dx<C<\infty.
\ee
Suppose also $\int^\pi_{-\pi} u^\theta (\cdot, z) dz = \int^\pi_{-\pi} u^z (\cdot, z) dz =0$.
Then $u=0$.
\end{theorem}

\begin{corollary}
\lab{codd}
Let $u$ be an axially symmetric solution to the problem
\be\label{e10.1}
\left\{
\begin{aligned}
&(u\cdot\nabla)u+\nabla p-\Delta u=0, \quad \text{in} \quad  \bR^2 \times [0, \pi],\\
&\nabla\cdot u=0,\\
&u^\theta |_{z=0, \pi}=0, \q  u^z |_{z=0, \pi}=0,\q  \p_z u^r |_{z=0, \pi}=0\\
&\lim\limits_{|x|\rightarrow \i}u=0,
\end{aligned}
\rt.
\ee
such that the Dirichlet integral satisfies the condition: for a constant $C$, and all $R \ge 1$,
\be
\int^\pi_{-\pi}\int_{R \le |x'| \le 2 R}|\nabla u(x)|^2 dx<C<\infty.\nn
\ee
Then $u=0$.
\end{corollary}

{\remark Clearly, if the Dirichlet integral is finite i.e. $\Vert \nabla u \Vert_{L^2(\bR^3 )} $ or $\Vert \nabla u \Vert_{L^2(\bR^2 \times [-\pi, \pi])}$ is finite, then  \eqref{dr2r1} or (\ref{dr2r}) is satisfied. With some extra work, then one can reach the vanishing result assuming the integral in (\ref{dr2r}) grows at certain power of $R$, namely
  $\int^\pi_{-\pi}\int_{|x'|\leq r}|\nabla u(x)|^2dx<(1+r)^\alpha$ for some suitable and positive $\alpha$.
 In a subsequent paper \cite{CPZZ:1}, the extra condition in Theorem \ref{thzp} that $u^z, u^\th$ have zero mean in the $z$ direction has been removed, under the stronger assumption that the Dirichlet integral is finite. Also the current method, being much different from that one, potentially allows application for flows with infinite Dirichlet energy such as Kolmogorov flows.}

Now we outline the proof of the above results briefly. We start with the observation that in a dyadic ball, after scaling, the axially symmetric Navier-Stokes equation resembles a 2 dimensional one. Then the $Brezis-Gallouet\ inequality$ introduced in \cite{BG:1} implies that a smooth vector field with finite Dirichlet energy is almost bounded. After returning to the original scale, one can show that $u$ is bounded by $C \big(\f{\ln r}{r}\big)^{1/2}$ for large $r$. Next by combining the equations for $w$ in \eqref{1.5}, $Brezis-Gallouet\ inequality$ and scaling technique, we will show that, with $z$ taken as $0$ for convenience,
\be\label{1.9}
|w^\th(r,0)|\leq Cr^{-1}(\ln r)^{1/2}\|(u^r,u^\th,u^z)\|^{1/2}_{L^\i([\f{3}{4}r,\f{5}{4}r]\times[-r,r])},
\ee
and
\be\label{1.10}
\begin{aligned}
|w^r(r,0)|+|w^z(r,0)|&\leq Cr^{-1}(\ln r)^{1/2}\|(u^r,u^z)\|^{1/2}_{L^\i([\f{3}{4}r,\f{5}{4}r]\times[-r,r])}\\
&\q\ +Cr^{-1/2}(\ln r)^{1/2}\|(\na u^r,\na u^z)\|^{1/2}_{L^\i([\f{3}{4}r,\f{5}{4}r]\times[-r,r])}.
\end{aligned}
\ee
The details are given in \eqref{3.8e} and \eqref{3.9ee}. Then using the decay of $u$ and \eqref{1.9}, we can deduce that the decay rate of $w^\th$ is $r^{-5/4}(\ln r)^{3/4}$.

In order to obtain decay of $w^r$ and $w^z$ from \eqref{1.10}, we need the decay of $\na u^r,\na u^z$ which can be connected with $w^\th$ by the $Biot-Savart\ law$  $-\Dl (u^re_r+u^ze_z)=\na\times(w^\th e_\th)$. Then $\na u^r,\na u^z$ can be written as integral representations of $w^\th$, $\int_{\bR^3} K(x,y)w^\th(y)dy$, where $K(x,y)$ are Calderon-Zygmund kernels. The decay relations between $\na u^r,\na u^z$ and $w^\th$ are shown in Lemma \ref{l3.2l}. At last a combination of \eqref{1.10}, decay of $u$ and $\na u^r,\na u^z$ imply the decay of $w^r,w^z$ in Theorem \ref{thfullwdecay}.

In the $z-$periodic case such that $D=\bR^2\times S^1$,  we will show that the decay rate of the velocity $u$ is $r^{-(\f{3}{2})^-}$ for large $r$ which implies that $u\equiv 0$ by the work \cite{CSTY:1} and \cite{KNSS:1}. Several steps are needed to get the decay of $u$. As a preparation, the Green's function $G$ on $\bR^2\times S^1$ for those functions whose integral on $S^1$ is zero will be introduced and a series of properties of $G$ will be displayed. The key point is that $G$ and its gradient have exponential decay near infinity and the radial
component of the velocity $u^r$ is the $z$ derivative of the angular stream function $L^\th$. This rapid decay
makes the difference between the periodic case and the full space case where the decay rate of $G$ is polynomial.

The details are done in a few steps. In step one: first we will use the $Biot-Savart\ law$ to get the representations of $u^r,u^z$ by $G$ and $w^\th$, displayed in \eqref{10.14} and \eqref{e4.36}, which indicate that $u^r,u^z$  decay in the same rate as $w^\th$ modulo a log term and an exponentially-decay  term
 \be\label{1.11}
\|(u^r,u^z)\|_{L^\i([\f{3}{4}r,\f{5}{4}r]\times[-\pi,\pi])}\leq C\ln r\|w^\th\|_{L^\i([\f{1}{2}r,\f{3}{2}r]\times[-\pi,\pi])}+Ce^{-\f{r}{64}}.
 \ee
 Then using $Brezis-Gallouet\ inequality$ and scaling technique, we will show that
\be\label{1.12}
\begin{aligned}
&\q\ \|(w^r,w^z)\|_{L^\i([\f{7}{8}r,\f{9}{8}r]\times[-\pi,\pi])}\\
&\leq C(\ln r)^{1/2}\|(u^r,u^z)\|_{L^\i([\f{3}{4}r,\f{5}{4}r]\times[-\pi,\pi])}.
\end{aligned}
\ee
Then a simple observation indicates that
\be\label{1.13}
 \|u^\th\|_{L^\i([\f{7}{8}r,\f{9}{8}r]\times[-\pi,\pi])}\leq \|w^r\|_{L^\i([\f{7}{8}r,\f{9}{8}r]\times[-\pi,\pi])}.
\ee
Using almost the same technique as the one for \eqref{1.9}, we can deduce
\be\label{1.14}
\|w^\th\|_{L^\i([\f{15}{16}r,\f{17}{16}r])}\leq Cr^{-1/2}(\ln r)^{1/2}\|(u^r,u^\th,u^z)\|^{1/2}_{L^\i([\f{7}{8}r,\f{9}{8}r]\times[-\pi,\pi])}.
\ee
The above four estimates allow us to improve the decay of $u$ and $w$ to $r^{-1^-}$ after  finitely many  iterations.

In step two, we differentiate the equations of $w$ to get the equations for $\na w$. By using the above estimates, we can prove that the decay rate of $|\na w|$ will be $r^{-(\f{3}{2})^-}$.

At last, we use the representations of $u^r,u^z$ again and another observation for $u^\th$ to show that the decay rate of $u$ is $r^{-(\f{3}{2})^-}$.

As an application of Theorem \ref{thzp},  to prove the vanishing of $u$ in \eqref{e10.1}, we perform even extension for $u^r$, the pressure $p$, and odd extension for $u^\th,u^z$  with respect to the $z$ variable. After such extensions, it is easy to see that $(u^r,u^\th,u^z)$ is a solution of \eqref{e10.2} satisfying the assumption in Theorem \ref{thzp}. Consequently, we can prove Corollary \ref{codd}.

Using the method described above, we can further prove a similar decay estimate for full 3 dimensional D solutions except for a small set of  angles.

\begin{theorem}
\lab{th3ddecay}
Let $u$ be a smooth  solution to the problem
\be
\left\{
\begin{aligned}
&(u\cdot\nabla)u+\nabla p-\Delta u=0, \quad \text{in} \quad  \bR^3,\\
&\nabla\cdot u=0,\\
&\lim\limits_{|x|\rightarrow \i}u=0,
\end{aligned}
\rt.
\ee
such that the Dirichlet integral satisfies the condition: for a constant $C$, and all $R \ge 1$,
\be\label{di<c}
\int^\i_{-\i}\int_{R \le |x'| \le 2 R}\big(|\nabla u(x)|^2 dx+|u(x)|^6\big) dx'dx_3<C<\infty.
\ee

 Let $\{ e_r, e_z, e_\theta\}$ be any given cylindrical system in which $x= (r, z, \theta) \in \bR^3$. Then for any $\delta \in (0, 1)$,  $r \in (R/2, 2 R)$, $ R \ge 2$ there exists  a measurable set $E  \equiv  E_R \subset [0, 2 \pi]$ with $|[0, 2 \pi]- E|<\delta$,  and a constant $C_0>0$, depending only on the constant $C$ in (\ref{di<c}) such that,
 \[
 |u(r, z, \th)| \le \frac{C_0}{\delta^{1/2}}  \left(\frac{\ln r}{ r} \right)^{1/2}, \qquad \forall \th \in E, \quad r \ge e, \, z \in \bR.
\]

\end{theorem}

\begin{remark} We should mention that
Robert Finn \cite{Fr:1961ACTA} (page229), already proved that for any vector field having a finite Dirichlet integral and tending to a constant vector $u_\i$ as $|x|\rightarrow \i$, then for any $\dl>0$, there exists a measurable set $E_\dl\subset\bS^2$, such that $|E_\dl|\leq \dl$ and
\be\label{rf1}
|u-u_\i|\leq \f{C}{\dl^{1/2}|x|^{1/2}}, \q \text{for}\ \forall x=|x|\o,\ \o\in\bS^2\backslash E_\dl.
\ee
This result and our result in the above theorem do not contain each other.  On one hand, we have a worse term $(\ln r )^{1/2}$. On the other hand, our decay estimate is uniform in all $z$ variable and the Dirichlet bounded integral is replaced by the more general \eqref{di<c}.
 The exception set is of two dimensional in \cite{Fr:1961ACTA} and our exception set is one dimensional. We wish to thank one referee for informing this.
\end{remark}

We conclude the introduction with a list of frequently used notations.  For $(x, t) \in \bR^3 \times \bR$ and $r>0$, we use $Q_r(x, t)$ to denote the standard parabolic cube $\{ (y, s) \q | \q |x-y|<r, 0<t-s<r^2 \}$.
The symbol $ ... \ls ...$ stands for $... \le C ...$ for a positive constant $C$. $B(x, r)$ denotes the ball of radius $r$ centered at $x$. $C$ with or without an index denotes a positive constant whose value may change from line to line.

\section{Mean value inequalities for velocity and decay in general domains}

In this section, we present mean value inequalities for the velocity for local smooth D-solutions of the Navier-Stokes equation in 3 dimensions. It is helpful in proving boundedness of local solutions without
explicit dependence on the pressure term. It is especially useful in proving a priori decay of D-solutions in the axially symmetric case with little restriction on the domain.  The result is an extension of those in \cite{O:1} and \cite{Z:2}  and its addendum, Lemma 5.1.

To start with, let us recall a basic inequality from \cite{O:1} Proposition 5, as modified in
the addendum of \cite{Z:2} since the former omitted the under-braced term in (\ref{2.2}) in the time dependent case. Since D solutions are independent of time, our proof is not affected.

\begin{lemma}
\lab{le2.1}
Let $\tu$ be a smooth solution of the Navier-Stokes equation \eqref{2.1} in $Q_{2r}(x,t)\subset\bR^3\times[0,\infty)$,
\be
\left\{
\begin{aligned}
&\p_t\tu+(\tu\cdot\nabla)\tu+\nabla p-\Delta\tu=0,\\
&\nabla\cdot\tu=0.
\end{aligned}
\rt. \label{2.1}
\ee

 Then there exists an absolute constant $\la$, independent of $r$ or $u$, such that
\be\label{2.2}
\begin{aligned}
|\tu(x,t)|&\leq \la \f{1}{r^5}\int_{Q_r(x,t)-Q_{r/2}(x,t)}|\tu(y,s)|dyds
 + \underbrace{\la  \f{1}{r^3}\int_{B_r(x)-B_{r/2}(x)}|\tu(y,t)|dy} \\
           &\q\ +\la\int_{Q_r(x,t)}K(x,t;y,s)|\tu(y,s)|^2dyds,
\end{aligned}
\ee
where $K(x,t;y,s)$ is defined as $K(x,t;y,s)=\big(|x-y|+\s{t-s}\big)^{-4}$.
\end{lemma}

Next we present a mean value inequality for a solution of \eqref{1.1}, which is not necessarily
axially symmetric.

\begin{proposition} [Mean value inequality]
\label{mean}
Let $u $ be a solution of \eqref{1.1}.  Then there exists a sufficiently small constant $\delta_0$ and another constant $C$ with the following property.   If
$B(x, 2 r_0) \subset D$ and $0<r_0 \le \delta_0  \Vert  u  \Vert^{-2}_{L^6(D)}$, then
\be
|u(x)|\leq \f{C}{r^3_0}\int_{B(x, r_0)}|u(y)|dy,  \nn\\
\ee
where $C$ is independent of $r_0$.
\end{proposition}

\q Before the proof of Proposition \ref{mean}, we need to recall, without proof,  a short lemma in \cite{Z:1}, concerning the kernel function $K$ defined in Lemma \ref{le2.1}.

\begin{lemma}
\lab{lez1}
For any locally integrable scalar function $u=u(y,s)$ and $t_1>t_2$, set
\be
\mathcal{K}(u;t_1,t_2)\triangleq \sup\limits_{x\in\bR^3}\int^{t_1}_{t_2}\int_{\bR^3}\big(K(x,t_1;y,s)+K(x,s;y,t_2)\big)|u(y,s)|dyds. \nn
\ee

Then, for $x,y,z\in\bR^3$ and $t_1>s>\tau>t_2$, we have
\be \label{3.1}
\begin{aligned}
&\int^{t_1}_{t_2}\int_{\bR^3}\big(|x-y|+\s{t_1-s}\big)^{-4}\big(|y-z|+\s{s-\tau}\big)^{-4}|u(y,s)|dyds \\
              & \leq C\mathcal{K}(u;t_1,t_2)K(x,t_1;z,\tau),
\end{aligned}
\ee
where $C$ is independent of $t_1, t_2$.
\end{lemma}

\noindent\textbf{Proof of Proposition \ref{mean}.}  For simplicity we replace $r_0$ by $r$ in the proof.  First we regard the solution $u=u(x)$ as a stationary solution of the Navier-Stokes equation
 (\ref{2.1}) in space-time domain $D \times[0,\infty)$. i.e. Set $\tu(x,t)=u(x)$. Then $\tu(x,t)$ is a solution to (\ref{2.1}) in $D \times [0, \infty)$. 
  Since the solution is independent of time, the two terms on the right hand side of line 1 in (\ref{2.2}) merge into one term. Here the constant $\la$ may be adjusted.

To simplify the notation, we use capital letters to denote points in space and time. Such as $X=(x,t), Y=(y,s), Z=(z,\tau)$. From Lemma 2.1, we have
\be
|\tu(X)|\leq\la\f{1}{r^5}\int_{Q_r(X)}|\tu(Y)|dY+\la\int_{Q_r(X)}K(X;Y)|\tu(Y)|^2dY, \label{3.2}
\ee

For $Y\in Q_r(X)$, we apply \eqref{3.2} in $Q_{r/2}(Y)$ and get
\be
|\tu(Y)|\leq \la\f{2^5}{r^5}\int_{Q_{r/2}(Y)}|\tu(Z)|dZ+\la\int_{Q_{r/2}(Y)}K(Y;Z)|\tu(Z)|^2dZ \label{3.3}
\ee
Inserting \eqref{3.3} into the second term of the right hand of \eqref{3.2}, we obtain
\bea
&|\tu(X)|&\leq\nn \\
        &&\la\f{1}{r^5}\int_{Q_r(X)}|\tu(Y)|dY+\la\int_{Q_r(X)}K(X;Y)|\tu(Y)|\la\f{2^5}{r^5}\int_{Q_{r/2}(Y)}|\tu(Z)|dZdY\nn  \\
        &&+\la\int_{Q_r(X)}K(X;Y)|\tu(Y)|\la\int_{Q_{r/2}(Y)}K(Y;Z)|\tu(Z)|^2dZdY \nn
\eea
Note that $Q_{r/2}(Y)\subset Q_{3r/2}(X)$ when $Y\in Q_r(X)$. Then we have
\be \label{3.4}
\begin{aligned}
|\tu(X)|&\leq
        \la\f{1}{r^5}\int_{Q_r(X)}|\tu(Y)|dY+\la^2\f{2^5}{r^5}\|\tu\|_{L^1(Q_{2r}(X))}\int_{Q_r(X)}K(X;Y)|\tu(Y)|dY  \\
        &+\la^2\int_{Q_{3r/2}(X)}\int_{Q_r(X)}K(X;Y)K(Y;Z)|\tu(Y)|dY|\tu(Z)|^2dZ
\end{aligned}
\ee
Applying Lemma \ref{lez1}, where we choose $t_1=t$ and $t_2=t-4r^2$ and $u=|\tu|$, we can get
\be
\int_{Q_r(X)}K(X;Y)K(Y;Z)|\tu(Y)|dY\leq c\mathcal{K}(|\tu|;t,t-4r^2)K(X;Z). \nn
\ee
Inserting this into \eqref{3.4}, we obtain
\be\lab{ux1<}
\begin{aligned}
|\tu(X)|&\leq \la\f{1}{r^5}\|\tu\|_{L^1(Q_{2r}(X))}+\la^2\f{2^5}{r^5}\|\tu\|_{L^1(Q_{2r}(X))}\mathcal{K}(|\tu|)  \\
        &+\la^2c\mathcal{K}(|\tu|)\int_{Q_{3r/2}(X)}K(X;Z)||\tu(Z)|^2dZ
        \end{aligned}
\ee
where $\mathcal{K}(|\tu|)=\mathcal{K}(|\tu|;t,t-4r^2)$.

For $\tu(Z)$, we use the inequality \eqref{3.2} in $Q_{r/4}(Z)$:
\be
|\tu(Z)|\leq\la\f{4^5}{r^5}\int_{Q_{r/4}(Z)}|\tu(W)|dW+\la\int_{Q_{r/4}(X)}K(Z;W)|\tu(W)|^2dW. \nn
\ee
Note that for $Z\in Q_{3r/2}(X)$, we have $Q_{r/4}(Z)\subset Q_{7r/4}(X)\subset Q_{2r}(X)$.
Hence, we can substitute the preceding inequality into the last term of (\ref{ux1<}) to
deduce
\bea
|\tu(X)|&\leq& \la\f{1}{r^5}\|\tu\|_{L^1(Q_{2r}(X))}+\la^2\f{2^5}{r^5}\|\tu\|_{L^1(Q_{2r}(X))}\mathcal{K}(|\tu|)\nn  \\
        &&+\la^3\f{4^5}{r^5}\|\tu\|_{L^1(Q_{2r}(X))}c\mathcal{K}(|\tu|)\int_{Q_{3r/2}(X)}K(X;Z)|\tu(Z)|dZ \nn \\
        &&+\la^3c\mathcal{K}(|\tu|)\int_{Q_{3r/2}(X)}K(X;Z)K(Z;W)|\tu(Z)|dZ\int_{Q_{r/4}(Z)}|\tu(W)|^2dW.\nn \\
        &\leq&\la\f{1}{r^5}\|\tu\|_{L^1(Q_{2r}(X))}+\la^2\f{2^5}{r^5}\|\tu\|_{L^1(Q_{2r}(X))}\mathcal{K}(|\tu|)\nn  \\
        &&+\la^3\f{4^5}{r^5}\|\tu\|_{L^1(Q_{2r}(X))}c\big(\mathcal{K}(|\tu|)\big)^2 \nn \\
        &&+\la^3\big(c\mathcal{K}(|\tu|)\big)^2\int_{Q_{7r/4}(X)}K(X;W)|\tu(W)|^2dW.\nn
\eea
Iterating this process and halving the size of cubic each time, by induction, we arrive at
\be
|\tu(X)|\leq \f{C}{r^5}\|\tu\|_{L^1(Q_{2r}(X))}\sum\limits^{\infty}_{i=1}\big(2^5c\mathcal{K}(|\tu|)\big)^i. \label{3.5}
\ee Here the constant in the sum $c$ may contain $\la$ as a factor.

Now we come to compute $\mathcal{K}(|\tu|)$.
\be
\begin{aligned}
\mathcal{K}(|\tu|)&=\mathcal{K}(|\tu|;t,t-4r^2) \\
        &=\sup\limits_{x\in\bR^3}\int^{t}_{t-4r^2}\int_{\bR^3}\big(K(x,t;y,s)+K(x,s;y,t-4r^2)\big)|\tu(y,s)|dyds \\
        &=\sup\limits_{x\in\bR^3}\int^{t}_{t-4r^2}\int_{\bR^3}\big(|x-y|+\s{t-s}\big)^{-4}|\tu(y,s)|dyds\\
        &\q\ +\sup\limits_{x\in\bR^3}\int^{t}_{t-4r^2}\int_{\bR^3}\big(|x-y|+\s{s-t+4r^2}\big)^{-4} |\tu(y,s)|dyds \\
        &\leq 2\sup\limits_{x\in\bR^3}\int^{4r^2}_0\Big(\int_{\bR^3}|u(y)|^6dy\Big)^{1/6}\Big(\int_{\bR^3}\big(|x-y|+\s{s}\big)^{-24/5}dy\Big)^{5/6}ds\\
        &\leq 2C\Big(\int_{\bR^3}|u(y)|^6dy\Big)^{1/6}\int^{4r^2}_0\Big(\int_{\bR^3}\big(|y|+\s{s}\big)^{-24/5}dy\Big)^{5/6}ds\\
        &\leq 2Cr^{1/2}\|u\|_{L^6(\bR^3)}.
\end{aligned}
\ee  In the above, if the domain of $u$ is $D$, then $\bR^3$ should be replaced by $D$.

We can choose a sufficiently small $\delta_0<1$,  such that when $r<\delta_0 \|u\|^{-2}_{L^6(\bR^3)}$, the following holds
 \[
 2^5c \mathcal{K}(|\tu|)< 1.
 \]
Then the series on the right hand side \eqref{3.5} converges, implying
\be
|\tu(X)|\leq \f{C}{r^5}\|\tu\|_{L^1(Q_{2r}(X))}.  \nn
\ee
Remember that $\tu(x,t)=u(x)$ actually is a function of $x$, independent of $t$. Then we have
\bea
|u(x)|&\leq& \f{C}{r^3}\int_{B(x, r)}|u(y)|dy.    \label{3.7}
\eea \qed

As an application, we prove an a priori decay estimate for general axially symmetric D-solutions. Notice that the domain can be any axially symmetric, noncompact one. When the domain is $\bR^3$, better estimate exists see \eqref{cjw} (c.f. \cite{CJ:1}, \cite{We:1}).

\begin{corollary} Let $D \subset \bR^3$ be a noncompact,connected and axially symmetric domain.
Suppose $B(x, 2) \subset D$ with $|x'|=r$. Then
there exists a constant $C$ independent of $u$ and $r$  such that
\[
|u(x)| \le \frac{C}{r^{1/6}} \Vert u \Vert^{5/3}_{L^6(D)}.
\]
\end{corollary}

\proof Fixing a radius $r_0= \min \{ 0.5 \delta_0 \|u\|^{-2}_{L^6(D)}, 1 \}$ such that Proposition \ref{mean} holds. Then
\be
|u(x)|\leq \f{C}{r^3_0}\int_{B_{r_0}(x)}|u(y)|dy \le \f{C}{r^{1/2}_0} \left(\int_{B_{r_0}(x)}|u(y)|^6 dy \right)^{1/6}.  \nn\\
\ee Since $u$ is axially symmetric, using the rotation method in \cite{LNZ:1}, we know that
\[
\int_{B_{r_0}(x)}|u(y)|^6 dy \le C \frac{r_0}{r} \int_D |u(y)|^6 dy.
\]Hence
\[
|u(x)| \le \f{C}{r^{1/2}_0} \frac{C r^{1/6}_0}{r^{1/6}} \| u\|_{L^6(D)} \le \frac{C}{r^{1/6}} \| u \|^{5/3}_{L^6(D)}.
\]\qed

\section{A priori decay of vorticity and velocity in the full space case}

The goal of the section is to give a proof Theorem \ref{thfullwdecay}, which will take a few steps.

We will frequently use the following notations.  If $a, b, c$ are scalars, then
$(a, b)$, $(a, b, c)$ denote, respectively, a 2-tuple and a 3-tuple, and $|(a, b)|=|a|+|b|$, etc. If $x=(x_1, x_2, x_3) \in \bR^3$, then $r=|x'|=\sqrt{x^2_1+x^2_2}$.

\subsection{Decay of $u^r, u^{\theta}, u^z$:
$|(u^r, u^{\theta}, u^z)|(x)\leq C \big(\f{\ln r}{r}\big)^{1/2}$}\ \\

As mentioned earlier, this bound was proven in \cite{CJ:1}, \cite{We:1} when the Dirichlet integral is finite. Here we give a very short proof, which shows that this type of decay is typical for smooth axially symmetric vector fields  and there is almost no use of equations.

Fixing $x_0 \in \bR^3$ such that $|x'_0|=  r_0$ is large. Without loss of generality, we can assume, in the cylindrical coordinates, that $x_0=(r_0, 0, 0)$, i.e. $z_0=0$, $\theta_0=0$.
Consider the scaled  solution
\[
\t {u} ( \t{x}) = r_0 u( r_0 \t{x})
\]which is also axially symmetric. Hence $\t {u}$ can be regarded as a two variable function  of the scaled variables $\t{r}, \t{z}$. Consider the two dimensional domain
\[
\t D =\{ (\t{r}, \t{z}) \, |,
1/2 \le \t{r} \le 2, \,  |\t{z}| \le 1 \}.
\]Then for $\t{u}=\t{u}(\t{r}, \t{z})$, we have $\t{u}(1, 0)=r_0 u(x_0)$.

To proceed,  we need the $Brezis-Gallouet$ inequality (\cite{BG:1}):

\begin{lemma}\label{l3.1}
Let $f\in H^2(\O)$ where $\O \subset \bR^2$. Then there exists a constant $C_{\O}$, depending only on $\O$, such that
\be\label{BGinq}
\|f\|_{L^\i(\O)}\leq C_{\O}\|f\|_{H^1(\O)}\ln^{1/2} \big(e+\f{\|\Dl f\|_{L^2(\O)}}{\|f\|_{H^1(\O)}}\big).
\ee
\end{lemma}

We mention that the original $Brezis-Gallouet$ inequality can be written in the form of
\[
\|f\|_{L^\i(\O)}\leq C_{\O}\|f\|_{H^1(\O)}\ln^{1/2} \big(e+\f{\| f\|_{H^2(\O)}}{\|f\|_{H^1(\O)}}\big).
\]However, by going through the proof on p678 of \cite{BG:1}, one can see that the norm $\| f\|_{H^2(\O)}$ in the log term can be replaced by $\|\Dl f\|_{L^2(\O)} + \| f\|_{L^2(\O)}$. Hence \eqref{BGinq} is valid.

 After a simple adjustment on the constant in \eqref{BGinq}, considering separately  the cases $\|f\|_{H^1(\O)} \le 1 $ and $\|f\|_{H^1(\O)} > 1 $,
 we  find that \eqref{BGinq} can be replaced by
\be\label{BGinq1}
\|f\|_{L^\i(\O)}\leq C_{\O}\Big(\|f\|_{H^1(\O)}+1\Big)\ln^{1/2} \big(e+\|\Dl f\|_{L^2(\O)}\big).
\ee Here is the detailed proof.  If $\|f\|_{H^1(\O)} > 1 $, then the above inequality clearly follows from (\ref{BGinq}).  If $\|f\|_{H^1(\O)} \le 1$, then write $a=\|f\|_{H^1(\O)}$ and $b
=  \|\Dl f\|_{L^2(\O)}$, we compute, from (\ref{BGinq}):
\[
a \ln ^{1/2} (e + (b/a)) \le a \ln^{1/2}(a e + b) + a \ln^{1/2} (1/a) \le C \ln^{1/2}(e+b).
\]which implies (\ref{BGinq1}).

Applying the Brezis-Gallouet inequality  \eqref{BGinq1} on $\t{D}$ , we can find an absolute constant $C$ such that
\[
\al
| \t{u}(1, 0) | &\le C \left[ \left( \int_{\t{D}} |\t{\nabla} \t{u} |^2 d\t{r}d\t{z} \right)^{1/2} +
\left(\int_{\t{D}} | \t{u} |^2 d\t{r}d\t{z}\right)^{1/2} +1  \right] \times \\
&\qquad \ln^{1/2}\left[ \left( \int_{\t{D}} |\t{\Delta} \t{u} |^2 d\t{r}d\t{z} \right)^{1/2} +
\left(\int_{\t{D}} | \t{u} |^2 d\t{r}d\t{z}\right)^{1/2} + e  \right],
\eal
\]where $\t \nabla = (\partial_{\t{r}}, \partial_{\t{z}})$ and
$\t \Delta = \partial^2_{\t{r}} + \partial^2_{\t{z}}$. By H\"older inequality and the assumption that
$1/2 \le \t{r} \le 2$, we see that
\be
\lab{tu10<}
\al
| \t{u}(1, 0) | &\le C \left[ \left( \int_{\t{D}} |\t{\nabla} \t{u} |^2 \t{r} d\t{r}d\t{z} \right)^{1/2} +
\left(\int_{\t{D}} | \t{u} |^6 \t{r} d\t{r}d\t{z}\right)^{1/6} +1  \right] \times \\
&\qquad \ln^{1/2}\left[ \left( \int_{\t{D}} |\t{\Delta} \t{u} |^2 \t{r} d\t{r}d\t{z} \right)^{1/2} +
\left(\int_{\t{D}} | \t{u} |^2 \t{r} d\t{r}d\t{z}\right)^{1/2} + e  \right].
\eal
\ee Now we can scale this inequality back to the original solution $u$ and variables $r=r_0 \t{r}$ and $z=r_0 \t{z}$ to get
\[
\al
r_0 |u(x_0)| &\le C  \left[ \sqrt{r_0} \left( \int_{D_0} |\nabla u |^2 r drdz \right)^{1/2} + \sqrt{r_0}
\left(\int_{D_0} |u|^6 r drdz \right)^{1/6} +1  \right] \times \\
&\qquad \ln^{1/2}\left[ r^{3/2}_0 \left( \int_{D_0} (|\p^2_r u|^2 + |\p^2_z u|^2) r drdz \right)^{1/2} + r^{-1/2}_0
\left(\int_{D_0} | u |^2 r drdz \right)^{1/2} + e  \right],
\eal
\]where
\[
 D_0 =\{ (r, z) \, |
r_0/2 \le r \le 2 r_0, \,  |z|  \le r_0 \}.
\]Since the solution has bounded $C^2$ norm, from assumption (\ref{dr2r1}), this proves the claimed decay of velocity.

\subsection{Decay of $w$: $|w^\th(x)|\leq C r^{-5/4} (\ln r)^{3/4}$,\quad $|(w^r, w^z)|(x)\leq C r^{-9/8}(\ln r)^{11/8}$}\ \\

 {\bf Step 1.} Proof of $|w^\th|\leq C r^{-5/4} (\ln r)^{3/4}$ and a weaker decay of $(w^r, w^z)$.

 Picking $x_0 \in \bR^3$ such that $|x'_0|=  \la$ is large. Without loss of generality, we can assume, in the cylindrical coordinates, that $x_0=(\la, 0, 0)$, i.e. $z_0=0$, $\theta_0=0$.

Consider the scaled solution:
\bea
&&\t{u}(\t{x})=\la u(\la \t{x})=\la u(x)\nn\\
&&\t{w}(\t{x})=\la^2w(\la \t{x})=\la^2w(x)\nn
\eea
where $\t{x}=\f{x}{\la}$.

First for simplification of notation, we drop the $``\sim"$ for a moment when computations take place under the scaled sense. Select the domains
\be
\mathcal{C}_1=\{(r,\th,z):\f{1}{2}<r<\f{3}{2},\ 0\leq \th\leq 2\pi,\ |z|\leq 1\},\nn
\ee
\be
\mathcal{C}_2=\{(r,\th,z):\f{3}{4}<r<\f{5}{4},\ 0\leq \th\leq 2\pi,\ |z|\leq \f{1}{2}\}.\nn
\ee

Let $\psi(y)$ be a cut-off function satisfying $\sup \psi(y)\subset \mathcal{C}_1$ and $\psi(y)=1$ for $y\in\mathcal{C}_2$ such that the gradient of $\psi$ is bounded. Now testing the vorticity equation \eqref{1.5} with $w^r\psi^2$, $w^\th\psi^2$ and $w^z\psi^2$ respectively, we have

\be
-\int_{\mathcal{C}_1} w^r\psi^2(\Dl-\f{1}{r^2})w^r dy=-\int_{\mathcal{C}_1}\left[ (u^r\p_r+u^z\p_z)w^r\cdot w^r\psi^2 -(w^r\p_r+w^z\p_z)u^r\cdot w^r\psi^2 \right]dy. \nn
\ee
\be
-\int_{\mathcal{C}_1} w^\th\psi^2(\Dl-\f{1}{r^2})w^\th dy=-\int_{\mathcal{C}_1}\left[(u^r\p_r+u^z\p_z)w^\th\cdot w^\th\psi^2-\f{u^r}{r}(w^\th)^2 \psi^2+2\f{w^r}{r}u^\th w^\th\psi^2 \right] dy. \nn
\ee

\be
-\int_{\mathcal{C}_1} w^z\psi^2\Dl w^z dy=-\int_{\mathcal{C}_1}\left[(u^r\p_r+u^z\p_z)w^z\cdot w^z\psi^2-(w^r\p_r+w^z\p_z)u^z\cdot w^z\psi^2 \right]dy. \nn
\ee

Then we have
\be\label{4.1}
\begin{aligned}
&\q\ \int_{\mathcal{C}_1} \Big(|\nabla(w^r\psi)|^2+\f{(w^r)^2\psi^2}{r^2}\Big) dy  \\
&=\int_{\mathcal{C}_1}\Big((w^r)^2|\nabla\psi|^2-\f{1}{2}\psi^2(u^r\p_r+u^z\p_z)(w^r)^2\\
&\qq\qq\qq\qq\qq+(w^r)^2\psi^2\p_ru^r+w^rw^z\psi^2\p_zu^r\Big )dy\\
&=\int_{\mathcal{C}_1}\Big((w^r)^2|\nabla\psi|^2+\f{1}{2}(w^r)^2(u^r\p_r+u^z\p_z)\psi^2\\
&\qq\qq\qq\qq\qq -\f{u^r}{r}(w^r\psi)^2-2u^rw^r\psi\p_r(w^r\psi)-u^r\p_z(w^r\psi w^z\psi)\Big )dy\\
&\leq C(1+\|(u^r,u^z)\|_{L^\i(\mathcal{C}_1)})\|w^r\|^2_{L^2(\mathcal{C}_1)}+\f{1}{4}\big(\|\nabla(w^r\psi)\|^2_{L^2(\mathcal{C}_1)}+\|\nabla(w^z\psi)\|^2_{L^2(\mathcal{C}_1)}\big)\\
&\qq\qq\qq\qq\qq+C\|u^r\|^2_{L^\i(\mathcal{C}_1)}\big(\|w^r\|^2_{L^2(\mathcal{C}_1)}+\|w^z\|^2_{L^2(\mathcal{C}_1)}\big)\\
&\leq C(1+\|(u^r,u^z)\|^2_{L^\i(\mathcal{C}_1)})\|(w^r,w^z)\|^2_{L^2(\mathcal{C}_1)}+\f{1}{4}\|\big(\nabla(w^r\psi),\nabla(w^z\psi)\big)\|^2_{L^2(\mathcal{C}_1)},
\end{aligned}
\ee
or
\be\label{e4.1e}
\begin{aligned}
&\q\ \int_{\mathcal{C}_1} \Big(|\nabla(w^r\psi)|^2+\f{(w^r)^2\psi^2}{r^2}\Big) dy  \\
&=\int_{\mathcal{C}_1}\Big((w^r)^2|\nabla\psi|^2-\f{1}{2}\psi^2(u^r\p_r+u^z\p_z)(w^r)^2+(w^r)^2\psi^2\p_ru^r+w^rw^z\psi^2\p_zu^r\Big )dy\\
&=\int_{\mathcal{C}_1}\Big((w^r)^2|\nabla\psi|^2+\f{1}{2}(w^r)^2(u^r\p_r+u^z\p_z)\psi^2+(w^r)^2\psi^2\p_ru^r+w^rw^z\psi^2\p_zu^r\Big )dy\\
&\leq C\big(1+\|(u^r,u^z)\|_{L^\i(\mathcal{C}_1)}+\|(\na u^r,\na u^z)\|_{L^\i(\mathcal{C}_1)}\big)\|(w^r,w^z)\|^2_{L^2(\mathcal{C}_1)}.
\end{aligned}
\ee
\be\label{4.2}
\begin{aligned}
&\q\ \int_{\mathcal{C}_1} \Big(|\nabla(w^\th\psi)|^2+\f{(w^\th)^2\psi^2}{r^2}\Big) dy  \\
&=\int_{\mathcal{C}_1}\Big((w^\th)^2|\nabla\psi|^2-
\f{1}{2}(u^r\p_r+u^z\p_z)(w^\th)^2\psi^2+\f{u^r}{r}(w^\th)^2\psi^2-2\f{w^r}{r}u^\th w^\th\psi^2\Big )dy\\
&\leq C(1+\|(u^\th, u^r)\|_{L^\i(\mathcal{C}_1)})\|(w^r,w^\th)\|^2_{L^2(\mathcal{C}_1)}-\f{1}{2}\int_{\mathcal{C}_1}(u^r\p_r+u^z\p_z)(w^\th)^2\psi^2dy\\
&\leq C(1+\|(u^\th, u^r)\|_{L^\i(\mathcal{C}_1)})\|(w^r,w^\th)\|^2_{L^2(\mathcal{C}_1)}+\f{1}{2}\int_{\mathcal{C}_1}(w^\th)^2(u^r\p_r+u^z\p_z)\psi^2dy\\
&\leq C\big(1+\|(u^r,u^\th,u^z)\|_{L^\i(\mathcal{C}_1)}\big)\|(w^r,w^\th)\|^2_{L^2(\mathcal{C}_1)}.
\end{aligned}
\ee
\be\label{4.3}
\begin{aligned}
&\q\ \int_{\mathcal{C}_1} |\nabla(w^z\psi)|^2 dy  \\
&=\int_{\mathcal{C}_1}\Big((w^z)^2|\nabla\psi|^2-\f{1}{2}\psi^2(u^r\p_r+u^z\p_z)(w^z)^2\\
&\qq\qq\qq\qq\qq+w^rw^z\psi^2\p_ru^z+(w^z\psi)^2\p_zu^z\Big )dy\\
&=\int_{\mathcal{C}_1}\Big((w^z)^2|\nabla\psi|^2+\f{1}{2}(w^z)^2(u^r\p_r+u^z\p_z)\psi^2\\
&\qq\qq\qq\qq\qq-2u^zw^z\psi\p_z(w^z\psi){-\f{u^z}{r}w^rw^z}\psi^2-u^z\p_r(w^r\psi w^z\psi)\Big )dy\\
&\leq C(1+\|(u^r,u^z)\|_{L^\i(\mathcal{C}_1)})\|w^z\|^2_{L^2(\mathcal{C}_1)}+\f{1}{4}\big(\|\nabla(w^r\psi)\|^2_{L^2(\mathcal{C}_1)}+
\|\nabla(w^z\psi)\|^2_{L^2(\mathcal{C}_1)}\big)\\
&\qq\qq\qq\qq\qq+C\|u^z\|^2_{L^\i(\mathcal{C}_1)}\big(\|w^r\|^2_{L^2(\mathcal{C}_1)}+\|w^z\|^2_{L^2(\mathcal{C}_1)}\big)\\
&\leq C(1+\|(u^r,u^z)\|^2_{L^\i(\mathcal{C}_1)})\|(w^r,w^z)\|^2_{L^2(\mathcal{C}_1)}+
\f{1}{4}\|\big(\nabla(w^r\psi),\nabla(w^z\psi)\big)\|^2_{L^2(\mathcal{C}_1)}.
\end{aligned}
\ee
or
\be\label{e4.3e}
\begin{aligned}
&\q\ \int_{\mathcal{C}_1} |\nabla(w^z\psi)|^2 dy  \\
&=\int_{\mathcal{C}_1}\Big((w^z)^2|\nabla\psi|^2-\f{1}{2}\psi^2(u^r\p_r+u^z\p_z)(w^z)^2+w^rw^z\psi^2\p_ru^z+(w^z\psi)^2\p_zu^z\Big )dy\\
&=\int_{\mathcal{C}_1}\Big((w^z)^2|\nabla\psi|^2+\f{1}{2}(w^z)^2(u^r\p_r+u^z\p_z)\psi^2 +w^rw^z\psi^2\p_ru^z+(w^z\psi)^2\p_zu^z)\Big )dy\\
&\leq C(1+\|(u^r,u^z)\|_{L^\i(\mathcal{C}_1)}+\|(\na u^r,\na u^z)\|_{L^\i(\mathcal{C}_1)})\|(w^r,w^z)\|^2_{L^2(\mathcal{C}_1)}.
\end{aligned}
\ee
From \eqref{4.1}, \eqref{e4.1e} and \eqref{4.3} \eqref{e4.3e}, we obtain
\be \label{4.4}
\begin{aligned}
&\q\ \|\big(\nabla w^r, \nabla w^z\big)\|^2_{L^2(\mathcal{C}_2)}\\
 & \leq  C(1+\|(u^r,u^z)\|^2_{L^\i(\mathcal{C}_1)})\|(w^r,w^z)\|^2_{L^2(\mathcal{C}_1)},
 \end{aligned}
\ee
or
\be \label{e4.4e}
\begin{aligned}
&\q\ \|\big(\nabla w^r, \nabla w^z\big)\|^2_{L^2(\mathcal{C}_2)}\\
&\leq  C(1+\|(u^r,u^z)\|_{L^\i(\mathcal{C}_1)}+\|(\na u^r,\na u^z)\|_{L^\i(\mathcal{C}_1)})\|(w^r,w^z)\|^2_{L^2(\mathcal{C}_1)}.
\end{aligned}
\ee
From \eqref{4.2}, we obtain
\be\label{4.5}
\begin{aligned}
&\q\ \|\nabla w^\th\|^2_{L^2(\mathcal{C}_2)} \\
&\leq C \big(1+\|(u^r,u^\th,u^z)\|_{L^\i(\mathcal{C}_1)}\big)\|(w^r,w^\th)\|^2_{L^2(\mathcal{C}_1)}.
\end{aligned}
\ee

Now we set
\be
\bar{\mathcal{C}}_2:=\{(r,z):\f{3}{4}<r<\f{5}{4},\ |z|\leq 1/2\}.\nn
\ee

Applying Lemma \ref{l3.1} and using \eqref{e4.1e} , we have
\be\label{3.6e}
\begin{aligned}
&\q\ \|w^\th\|_{L^\i(\bar{\mathcal{C}}_2)} \\
&\ls \|w^\th\|_{H^1(\bar{\mathcal{C}}_2)}\ln^{1/2} \big(e+\f{\|\Dl w^\th\|_{L^2({\mathcal{C}}_2)}}{\|w^\th\|_{H^1(\bar{\mathcal{C}}_2)}}\big)\\
&\ls \big(1+\|w^\th\|_{H^1(\bar{\mathcal{C}}_2)}\big)\ln^{1/2} \big(e+\|\Dl w^\th\|_{L^2(\bar{\mathcal{C}}_2)}\big)\\
&\ls \Big(1+(1+\|(u^r,u^\th,u^z)\|^{1/2}_{L^\i(\mathcal{C}_1)})\|(w^r,w^\th)\|_{L^2({\mathcal{C}}_1)}\Big)
\ln^{1/2} \big(e+\|\Dl w^\th\|_{L^2({\mathcal{C}}_1)}\big);
\end{aligned}
\ee
Similarly applying Lemma \ref{l3.1} and using \eqref{4.4}, we find
\be\label{3.7e}
\begin{aligned}
&\q\ \|(w^r,w^z)\|_{L^\i(\bar{\mathcal{C}}_2)}  \\
&\ls \|(w^r,w^z)\|_{H^1(\bar{\mathcal{C}}_2)}\ln^{1/2} \big(e+\f{\|(\Dl w^r,\Dl w^z)\|_{L^2(\bar{\mathcal{C}}_2)}}{\|(w^r, w^z)\|_{H^1(\bar{\mathcal{C}}_2)}}\big)\\
&\ls \big(1+\|(w^r,w^z)\|_{H^1(\bar{\mathcal{C}}_2)}\big)\ln^{1/2} \big(e+\|(\Dl w^r,\Dl w^z)\|_{L^2(\bar{\mathcal{C}}_2)}\big)\\
&\ls \Big(1+(1+\|(u^r,u^z)\|_{L^\i(\mathcal{C}_1)})\|(w^r,w^z)\|_{L^2({\mathcal{C}}_1)}\Big)\\
&\q\ \times\ln^{1/2} \big(e+\|(\Dl w^r,\Dl w^z)\|_{L^2({\mathcal{C}}_1)}\big).
\end{aligned}
\ee
Similarly applying Lemma \ref{l3.1} and using \eqref{e4.4e}, we find
\be\label{3.7e1}
\begin{aligned}
&\q\ \|(w^r,w^z)\|_{L^\i(\bar{\mathcal{C}}_2)}  \\
&\ls \|(w^r,w^z)\|_{H^1(\bar{\mathcal{C}}_2)}\ln^{1/2} \big(e+\f{\|(\Dl w^r,\Dl w^z)\|_{L^2(\bar{\mathcal{C}}_2)}}{\|(w^r, w^z)\|_{H^1(\bar{\mathcal{C}}_2)}}\big)\\
&\ls \big(1+\|(w^r,w^z)\|_{H^1(\bar{\mathcal{C}}_2)}\big)\ln^{1/2} \big(e+\|(\Dl w^r,\Dl w^z)\|_{L^2(\bar{\mathcal{C}}_2)}\big)\\
&\ls\Big(1+(1+\|(u^r,u^z)\|^{1/2}_{L^\i(\mathcal{C}_1)}+\|(\na u^r,\na u^z)\|^{1/2}_{L^\i(\mathcal{C}_1)})\|(w^r,w^z)\|_{L^2({\mathcal{C}}_1)}\Big)\\
&\q\ \times\ln^{1/2} \big(e+\|(\Dl w^r,\Dl w^z)\|_{L^2({\mathcal{C}}_1)}\big).
\end{aligned}
\ee

Now scaling back, to the domains
\be
\mathcal{C}_{1, \la}=\{(r,\th,z):\f{\la}{2}<r<\f{3\la}{2},\ 0\leq \th\leq 2\pi,\ |z|\leq \la\},\nn
\ee,
\be
\mathcal{C}_{2, \la}=\{(r,\th,z):\f{3\la}{4}<r<\f{5\la}{4},\ 0\leq \th\leq 2\pi,\ |z|\leq \f{\la}{2}\},\nn
\ee

 we have

\be\label{3.8e}
\begin{aligned}
&\q\ \la^2\|w^\th\|_{L^\i(\mathcal{C}_{2,\la})}\\
&\ls\Big(1+\big(1+\la^{1/2}\|(u^r,u^\th,u^z)\|^{1/2}_{L^\i(\mathcal{C}_{1,\la})}\big)\la^{1/2}\|(w^r,w^\th)\|_{L^2(\mathcal{C}_{1,\la})}\Big)\times\\
&\q\ \ln^{1/2} \big(\la^{5/2}\|\Dl w^\th\|_{L^2(\mathcal{C}_{2,\la})}+e\big)\\
&\ls\la^{1/2}\Big(1+\big(1+\la^{1/2}\|(u^r,u^\th,u^z)\|^{1/2}_{L^\i(\mathcal{C}_{1,\la})}\big)\|(w^r,w^\th)\|_{L^2(\mathcal{C}_{1,\la})}\Big)\times\\
&\q\ \ln^{1/2} \big(\la^{5/2}\|\Dl w^\th\|_{L^2(\mathcal{C}_{2,\la})}+e\big)
\end{aligned}
\ee
and
\be\label{3.9e}
\begin{aligned}
&\q\ \la^2\|(w^r,w^z)\|_{L^\i(\mathcal{C}_{2,\la})}\\
&\leq C\la^{1/2}\Big(1+\big(1+\la\|(u^r,u^z)\|_{L^\i(\mathcal{C}_{1,\la})}
\big)\|(w^r,w^z)\|_{L^2(\mathcal{C}_{1,\la})}\Big)\\
&\q\ \times\ln^{1/2} \big(\la^{5/2}\|(\Dl w^r,\Dl w^z)\|_{L^2(\mathcal{C}_{2,\la})}+e\big).
\end{aligned}
\ee
\be\label{3.9ee}
\begin{aligned}
&\q\ \la^2\|(w^r,w^z)\|_{L^\i(\mathcal{C}_{2,\la})}\\
&\leq C\la^{1/2}\Big(1+\big(1+\la^{1/2}\|(u^r,u^z)\|^{1/2}_{L^\i(\mathcal{C}_{1,\la})}\\
&\qq\qq\qq+\la\|(\na u^r,\na u^z)\|^{1/2}_{L^\i(\mathcal{C}_{1,\la})}
\big)\|(w^r,w^z)\|_{L^2(\mathcal{C}_{1,\la})}\Big)\\
&\q\ \times\ln^{1/2} \big(\la^{5/2}\|(\Dl w^r,\Dl w^z)\|_{L^2(\mathcal{C}_{2,\la})}+e\big).
\end{aligned}
\ee

 We mention here that, according to  scaling and routine energy estimates, $\Delta w$ can only grow as a polynomial order at the far field. Thus we need not to calculate the exact order, since $\left\|\Delta w\right\|_{L^2(\mathcal{C}_{2,\lambda})}$ appears in a "$\log$" function.

From \eqref{3.8e}, \eqref{3.9e}, assumption (\ref{dr2r1}) and by the a priori bound on $u$, we have
\be\label{3.18r}
\|w^\th\|_{L^\i(\mathcal{C}_{2,\la})}\leq C\la^{-5/4}(\ln\la)^{3/4},\q \|(w^r,w^z)\|_{L^\i(\mathcal{C}_{2,\la})}\leq C\la^{-1}\ln\la
\ee
In addition, carrying out a similar estimate on the unit ball centered at $x \in \mathcal{C}_{2,\la}$, it is routine to show that
\be
\|\na w^\th\|_{L^\i(\mathcal{C}_{2,\la})}\leq C\la^{-5/4}(\ln\la)^{3/4}.\nn
\ee

\medskip

{\bf Step 2.}  Proof of $|(w^r, w^z)|\leq C r^{-9/8}(\ln r)^{11/8}$.

In order to get more decay estimates of $w^r,w^z$, we will use \eqref{3.9ee} and further $L^\i$ estimates of $\na u^r,\na u^z$. Note that
 according to Majda and Betozzi \cite{MB:2002CAMBRIDGE} (page 77)
\be\label{e3.16}
\na(u^re_r+u^ze_z)=\bC\cdot w^\th e_\th+\bK\ast(w^\th e_\th),
\ee
where $\bC$ is a constant matrix and $\bK$ is a Calderon-Zygmund kernel.

%which means that
%\be\label{e3.16}
%\na u^r=\int K_1(x,y)w_\th(y)dy,\q \na u^z=\int K_2(x,y)w_\th(y)dy,
%\ee
%where $K_1(x,y)$ and $K_2(x,y)$ are Calderon-Zygmund kernels.
It is also easy to check that the $y$ integral of kernels $\bK(x,y)$ on any balls centered at $x$ is $0$.

Next we present a lemma to describe a property of the Calderon-Zygmund kernel when it acts on some axially symmetric functions.

\begin{lemma}\label{l3.2l}
Assume that $K(x,y)$ be a Calderon-Zygmund kernel and $f$ is a smooth axisymmetric function satisfying, for $x=(x', z)\in \bR^3$
\be
|f(x)|+|\na f(x)|\leq \f{\ln^b(e+|x'|)}{(1+|x'|)^a}\q \text{for}\q 0<a<2,\ b>0.\nn
\ee
Define $Tf(x):=\int K(x,y)f(y)dy$.  Then there exists a constant $c_0$ such that
\be\label{ee3.17}
|Tf(x)|\leq c_0\f{\ln^{b+1}(e+|x'|)}{(1+|x'|)^a}.
\ee
\end{lemma}
\pf Note that for the Calderon-Zygmund kernel, we have
\be
P.V. \int_{|x-y|\leq 1}K(x,y)dy=0. \nn
\ee
So we decompose $Tf(x)$ as follows
\bea
Tf(x)&=&\int_{\{|x-y|\leq 1\}}K(x,y)(f(y)-f(x))dy\nn\\
&&+\int_{\{|x-y|\geq 1\}\cap \{|x'-y'|\geq 1/2\}}K(x,y)f(y)dy+\int_{\{|x-y|\geq 1\}\cap \{|x'-y'|\leq 1/2\}}K(x,y)f(y)dy\nn\\
&:=&\sum\limits^3_1 I_i.\nn
\eea
Then using mean-value inequality, we have
\be\label{e3.18e}
\begin{aligned}
|I_1|&\ls\int_{\{|x-y|\leq 1\}}|K(x,y)|{|x-y|}|\na f(\xi)|dy\q \ |\xi| \in(\min\{|x|,|y|\},\max\{|x|,|y|\})\\
&\ls\int_{\{|x-y|\leq 1\}}\f{1}{|x-y|^3}|x-y|\f{\ln^b(e+|\xi'|)}{(1+|\xi'|)^a}dy\\
&\ls \f{\ln^b(e+|x'|)}{(1+|x'|)^a}\int_{\{|x-y|\leq 1\}}\f{1}{|x-y|^2}dy\\
&\ls\f{\ln^b(e+|x'|)}{(1+|x'|)^a}.
\end{aligned}
\ee
For the term $I_2$, we have
\bea
I_2&\ls&\int_{\substack{\{|x-y|\geq 1\}\\\{|x'-y'|\geq 1/2\}}} \f{\ln^b(e+|y'|)}{(1+|y'|)^a} \int^{+\i}_{-\i}\f{1}{(|x'-y'|+|x_3-y_3|)^3}dy_3 dy'\nn\\
&\ls&\int_{\{|x'-y'|\geq 1/2\}} \f{1}{|x'-y'|^2}\f{\ln^b(e+|y'|)}{(1+|y'|)^a}dy'\nn\\
&=&\bigg(\int_{\substack{\{|x'-y'|\geq 1/2\}\\\{|y'|\geq 2|x'|\}}}+\int_{\substack{\{|x'-y'|\geq 1/2\}\\\{1/2|x'|\leq|y'|\leq 2|x'|\}}}+\int_{\substack{\{|x'-y'|\geq 1/2\}\\\{|y'|\leq 1/2|x'|\}}}\bigg) \f{1}{|x'-y'|^2}\f{\ln^b(e+|y'|)}{(1+|y'|)^a}dy'\nn\\
&:=&\sum^3_i I_{2,i}.\nn
\eea
Now we estimates $I_{2,i}$ term by term
\bea
I_{2,1}&\ls& \int_{\{|y'|\geq \max\{1/3,2|x'|\}\}}\f{\ln^b(e+|y'|)}{(1+|y'|)^{2+a}}dy'\nn\\
&\ls&\int_{\{|y'|\geq \max\{1/3,2|x'|\}\}}\f{\ln^b(e+|y'|)}{(1+|y'|)^{1+a}}d|y'|\nn\\
&\ls&\f{\ln^b(e+|x'|)}{(1+|x'|)^a}.\nn
\eea
 \bea
I_{2,2}&\ls& \f{\ln^b(e+|x'|)}{(1+|x'|)^a}\int_{\substack{\{|x'-y'|\geq 1/2\}\\\{1/2|x'|\leq|y'|\leq 2|x'|\}}}\f{1}{|x'-y'|^2}dy'\nn\\
&\ls&\f{\ln^b(e+|x'|)}{(1+|x'|)^a}\int_{\{1/2 |x'| \leq s\leq 3|x'|\}}\f{1}{s}ds\nn\\
&\ls&\f{\ln^{b+1}(e+|x'|)}{(1+|x'|)^a}.\nn
\eea
\bea
I_{2,3}&\ls& \f{\ln^b(e+|x'|)}{(1+|x'|)^2}\int_{\substack{\{|x'-y'|\geq 1/2\}\\\{|y'|\leq1/2|x'|\}}}\f{1}{(1+|y'|)^{a}}dy'\nn\\
&\ls& \f{\ln^b(e+|x'|)}{(1+|x'|)^2}\int_{\{|y'|\leq 1/2|x'|\}}\f{1}{(1+|y'|)^{a-1}}d|y'|\nn\\
&\ls&\f{\ln^{b}(e+|x'|)}{(1+|x'|)^a}.\nn
\eea
The above estimates indicate that
\be\label{e3.19e}
I_2\ls \f{\ln^{b+1}(e+|x'|)}{(1+|x'|)^a}.
\ee
For the term $I_3$, we have
\be\label{e3.20e}
\begin{aligned}
I_3&\ls\int_{\substack{\{|x-y|\geq 1\}\\\{|x'-y'|\leq 1/2\}}}\f{1}{|x-y|^3}\f{\ln^{b}(e+|y'|)}{(1+|y'|)^a}dy\\
&\ls\int_{\substack{\{|x_3-y_3|\geq 1/2\}\\\{|x'-y'|\leq 1/2\}}}\f{1}{|x-y|^3}\f{\ln^{b}(e+|y'|)}{(1+|y'|)^a}dy\\
&\ls\f{\ln^{b}(e+|x'|)}{(1+|x'|)^a}\int_{\{|x'-y'|\leq 1/2\}}\int_{\{|x_3-y_3|\geq 1/2\}}\f{1}{(|x'-y'|+|x_3-y_3|)^3}dy_3dy'\\
&\ls\f{\ln^{b}(e+|x'|)}{(1+|x'|)^a}.
\end{aligned}
\ee

Combining \eqref{e3.18e},\eqref{e3.19e} and \eqref{e3.20e}, we get \eqref{ee3.17}, finishing the proof of the lemma.\ef

Now applying Lemma \ref{l3.2l} to \eqref{e3.16} and noting the decay rate of $w^\th$ in \eqref{3.18r}, we can get
\be\label{e3.21}
|\na u^r|+ |\na u^z|\ls r^{-5/4}(\ln r)^{7/4} \q \text{for\ large}\ r.
\ee

Now we go back to \eqref{3.9ee}, we can get more decay on $w^r,w^z$
\be\label{e3.22e}
\begin{aligned}
|w^r|+|w^z|&\ls r^{-3/2}(r^{1/2}r^{-1/4}(\ln r)^{1/4}+r(r^{-5/4}(\ln r)^{7/4})^{1/2})(\ln r)^{1/2}\\
&\ls r^{-9/8}(\ln r)^{11/8} \q \text{for\ large}\ r.
\end{aligned}
\ee
 This proves Theorem \ref{thfullwdecay}. \qed

\subsection{Proof of Theorem \ref{th3ddecay}, decay of full 3d solutions except for a small set}
$$
$$

Now $u$ is a full 3 d solution.
Fixing $x_0 \in \bR^3$ such that $|x'_0|=  r_0$ is large. Without loss of generality, we can assume, in the cylindrical coordinates, that $x_0=(r_0, 0, \th)$.
Consider the scaled  solution
\[
\t {u} ( \t{x}) = r_0 u( r_0 \t{x}).
\] Fixing $\th$, then $\t {u}$ can be regarded as a two variable function of the scaled variables $\t{r}, \t{z}$. Consider the two dimensional domain
\[
\t D =\{ (\t{r}, \t{z}) \, |,
1/2 \le \t{r} \le 2, \,  |\t{z}| \le 1 \}.
\]Then for $\t{u}=\t{u}(\t{r}, \t{z}, \th)$, we have $\t{u}(1, 0, \th)=r_0 u(x_0)$.

Just as in Section 3.1, the inequality (\ref{tu10<}) continues to hold, namely,
\be
\lab{3tu10<}
\al
| \t{u}(1, 0, \th) | &\le C_1 \left[ \left( \int_{\t{D}} |\t{\nabla} \t{u} |^2(\t r, \t z, \th) \t{r} d\t{r}d\t{z} \right)^{1/2} +
\left(\int_{\t{D}} | \t{u} |^6(\t r, \t z, \th) \t{r} d\t{r}d\t{z}\right)^{1/6} +1  \right] \times \\
&\qquad \ln^{1/2}\left[ \left( \int_{\t{D}} |\t{\Delta} \t{u} |^2(\t r, \t z, \th) \t{r} d\t{r}d\t{z} \right)^{1/2} +
\left(\int_{\t{D}} | \t{u} |^2(\t r, \t z, \th) \t{r} d\t{r}d\t{z}\right)^{1/2} + e  \right].
\eal
\ee Here $\t \nabla =(\p_{\t r}, \p_{\t z})$. Note that there is no integration in the $\th$ direction. Now we can scale this inequality back to the original solution $u$ and variables $r=r_0 \t{r}$ and $z=r_0 \t{z}$ to get
\[
\al
r_0 |u(r_0, 0, \th)| &\le C_1  \left[ \sqrt{r_0} \left( \int_{D_0} |\nabla u |^2 (r, z, \th)r drdz \right)^{1/2} + \sqrt{r_0}
\left(\int_{D_0} |u|^6(r, z, \th) r drdz \right)^{1/6} +1  \right] \times \\
&\qquad \ln^{1/2}\left[ r^{3/2}_0 \left( \int_{D_0} (|\p^2_r u|^2 + |\p^2_z u|^2) r drdz \right)^{1/2} + r^{-1/2}_0
\left(\int_{D_0} | u |^2 r drdz \right)^{1/2} + e  \right],
\eal
\]where
\[
 D_0 =\{ (r, z) \, |
r_0/2 \le r \le 2 r_0, \,  |z|  \le r_0 \}.
\]Since the solution has bounded $C^2$ norm, this implies
\be
\lab{r0ur0<}
\al
&r_0 |u(r_0, 0, \th)| \\
&\le C_2  \left[ \sqrt{r_0} \left( \int_{D_0} |\nabla u |^2 (r, z, \th)r drdz \right)^{1/2} + \sqrt{r_0}
\left(\int_{D_0} |u|^6(r, z, \th) r drdz \right)^{1/6} +1  \right]  \ln^{1/2} r_0\\
&\le C_2  \sqrt{r_0} \left( \int^{r_0}_{-r_0} \int^{2 r_0}_{r_0/2} (|\nabla u |^2 + |u|^6)(r, z, \th)r drdz \right)^{1/2}   \ln^{1/2} r_0.
\eal
\ee Recall our assumption
\[
\int^{2 \pi}_0 \int^\infty_{-\infty}\int^{2r_0}_{r_0/2}
(
|\nabla u |^2 + |u|^6 )(r, z, \th)r drdzd\th \le C.
\]Hence, for any $\delta \in (0, 1)$, there exists a measurable set $E$ with $|[0, 2 \pi]-E|<\delta$ such that
\[
\int^\infty_{-\infty} \int^{2r_0}_{r_0/2} (|\nabla u |^2 + |u|^6) (r, z, \th)r drdz \le \frac{ C^2}{\delta}
\]for all $\th \in E$. Plugging this into (\ref{r0ur0<}), we deduce
\[
 |u(r_0, 0, \th)| \le C \left( \frac{ \ln r_0}{\delta r_0} \right)^{1/2}
\]for all $\th \in E$. This proves the theorem. \qed

\section{z-periodic D-solutions of axisymmetric Navier-Stokes equations with the assumption $\int^\pi_{-\pi}u^\th dz=\int^\pi_{-\pi}u^z dz=0$}\label{s4}

This section is divided into three subsections. The following notations are frequently used from now on.
For $x=(x_1, x_2, x_3) \in \bR^3$, we write $x=(x', x_3)$ or $x=(x', z)$, and for
$y=(y_1, y_2, y_3) \in \bR^3$, we write $y=(y', y_3)$.

\subsection{The Green function on $\bR^2\times S^1$ for functions whose integral on $S^1$ is zero}

\begin{lemma} Let $\tilde \Gamma_1$ and $\Gamma_2$ be the standard heat kernel on $S^1$ and $\bR^2$ respectively.

(a). Then the function
\be\label{4.1e}
G(x, y) = \int^\infty_0 (\tilde \Gamma_1-\frac{1}{2 \pi})(x_3, y_3, t) \Gamma(x', y', t) dt
\ee
is well defined except when $x=y$.

(b). Let $u$ and $f$ be smooth, bounded functions on $\bR^2\times S^1$ such that
\[
\Delta u = - f.
\]Suppose $\int_{S^1} f(x', x_3) dx_3=0$ for all $x' \in \bR^2$ and $f$ is compactly supported. Then
\[
u(x) = \int_{\bR^2\times S^1} G(x, y) f(y) dy + C.
\] Thus $G$ is the Green's function for those $f$.
\end{lemma}

\proof

\qq Let $S^1=[-\pi,\pi]$. Define $\t{\G}_1(t;x_3,y_3):=\f{1}{2\pi}\big(1+2\sum\limits_{m\in\bZ^+}e^{-m^2t}\cos(m(x_3-y_3))\big)$, where $(x_3,y_3)\in [-\pi,\pi]^2$. Then it is known  that $\t{\G}_1(t;x_3,y_3)$ is the heat kernel on $S^1$, which means that
\be
\lt\{
\begin{aligned}
&\Dl_{x_3}\t{\G}_1-\p_t\t{\G}_1=0,\\
&\t{\G}_1|_{t=0}=\dl(x_3-y_3).
\end{aligned}
\rt.\nn
\ee  For completeness we give a proof.
It is easy to check that $\t{\G}_1$ satisfies the equation, we only check that $\t{\G}_1|_{t=0}=\dl(x_3-y_3)$. For any test function $\varphi(y_3)\in C^\i_0(\bR)$,
\bea
&&\q\lim\limits_{t\rightarrow 0}\langle \t{\G}_1(t;x_3,y_3),\varphi(y_3)\rangle-\varphi(x_3)\nn\\
&&=\f{1}{2\pi} \lim_{t \to 0^+} \int^\pi_{-\pi}\big((1+2\sum\limits_{m\in\bZ^+}\cos(m(x_3-y_3)) e^{- m^2 t} \big)\varphi(y_3)dy_3-\f{1}{2\pi}\int^\pi_{-\pi}\varphi(x_3)dy_3\nn\\
&&=\f{1}{2\pi}\lim\limits_{n\rightarrow +\i}\int^\pi_{-\pi}\big((1+2\sum\limits^n_{m=1}\cos(m(x_3-y_3))\big)\varphi(y_3)dy_3-\f{1}{2\pi}\int^\pi_{-\pi}\varphi(x_3)dy_3\nn\\
&&=\f{1}{2\pi}\lim\limits_{n\rightarrow +\i}\int^\pi_{-\pi}\f{[\sin(n+1/2)(x_3-y_3)]}{\sin\f{x_3-y_3}{2}}\varphi(y_3)dy_3-\f{1}{2\pi}\int^\pi_{-\pi}\varphi(x_3)dy_3\nn\\
&&=\f{1}{2\pi}\lim\limits_{n\rightarrow +\i}\int^\pi_{-\pi}\f{\sin(n+1/2)(x_3-y_3)}{\sin\f{x_3-y_3}{2}}\big(\varphi(y_3)-\varphi(x_3)\big)dy_3\nn\\
&&=\f{1}{2\pi}\lim\limits_{n\rightarrow +\i}\int^\pi_{-\pi}\f{\varphi(y_3)-\varphi(x_3)}{\sin\f{x_3-y_3}{2}}\sin(n+1/2)(x_3-y_3)dy_3.\nn
\eea
Due to the fact
\be
\lim\limits_{x_3\rightarrow y_3}\f{\varphi(y_3)-\varphi(x_3)}{\sin\f{x_3-y_3}{2}}=2\varphi'(y_3).\nn
\ee
 $\f{\varphi(y_3)-\varphi(x_3)}{\sin\f{x_3-y_3}{2}}$ is integrable on $S^1$. By using Riemann-Lebesgue Lemma, we can get
\be
\lim\limits_{t\rightarrow 0}\langle \t{\G}(t;x_3,y_3),\varphi(y_3)\rangle-\varphi(x_3)=0.\nn
\ee

Now we define $\G_1(t;x_3,y_3)=\t{\G}_1(t;x_3,y_3)-\f{1}{2\pi}$. Then $\G_1(t;x_3,y_3)$ is still a heat kernel for the test function $\varphi(y_3)$ with $\int^\pi_{-\pi}\varphi(y_3)dy_3=0$ due to the fact
\bea
&&\lim\limits_{t\rightarrow 0}\langle {\G}_1(t;x_3,y_3),\varphi(y_3)\rangle-\varphi(x_3)\nn\\
&&=\lim\limits_{t\rightarrow 0}\langle \t{\G}_1(t;x_3,y_3),\varphi(y_3)\rangle-\varphi(x_3)-\f{1}{2\pi}\int^\pi_{-\pi}\varphi(y_3)dy_3\nn\\
&&=0.\nn
\eea

Next we construct the heat kernel on $\bR^2\times[-\pi,\pi]$ for bounded smooth functions $\varphi =\varphi(y)$ such that
\be\label{ee10.4}
\int^\pi_{-\pi}\varphi(y',y_3)dy_3=0.
\ee

Set $\G_2(t;x',y')=(4\pi t)^{-1}e^{-\f{|x'-y'|^2}{4t}}$ be the heat kernel on $\bR^2$.  Define
\be\label{ee4.2}
\G(t;x,y)=\G_2(t;x',y')\G_1(t;x_3,y_3).
\ee
 Then we claim that
\be
\lt\{
\begin{aligned}
&\Dl_x \G-\p_t \G=0\\
&\G|_{t=0}=\dl(x,y)
\end{aligned}
\rt.
\ee
with  bounded, smooth test functions satisfying \eqref{ee10.4}. Here the second line means $\lim_{t \to 0^+}
\Gamma * \varphi =\varphi$.

First it is clear that
\bea
&&\q\Dl_x \G-\p_t \G\nn\\
&&=\G_1\Dl_2\G_2+\G_2\Dl_1\G_1-\G_1\p_t\G_2-\G_2\p_t\G_1\nn\\
&&=0.\nn
\eea
Next for any bounded, smooth test function $\varphi(y)$ satisfying \eqref{ee10.4}, we have
\bea
&&\q\lim\limits_{t\rightarrow 0}\int_{\bR^2\times S^1}\G(t;x,y)\varphi(y)dy \nn\\
&&=\lim\limits_{t\rightarrow 0}\int_{\bR}\G_2(t;x',y')\int_{S^1}\G_1(t;x_3,y_3)\varphi(y',y_3)dy_3dy'\nn\\
&&=\lim\limits_{t\rightarrow 0}\int_{\bR}\G_2(t;x',y')\int_{S^1}\tilde \G_1(t;x_3,y_3)\varphi(y',y_3)dy_3dy'\nn\\
&&=\lim\limits_{t\rightarrow 0}\int_{\bR}\G_2(t;x',y')\varphi(y',x_3)dy'
+ \lim\limits_{t\rightarrow 0}\int_{\bR}\G_2(t;x',y') [-\varphi(y', x_3)\nn\\
&&\qq\qq+ \int_{S^1}\tilde \G_1(t;x_3,y_3)\varphi(y',y_3)dy_3 ] dy' \nn\\
&&=\varphi(x',x_3)\nn\\
&&=\varphi(x)\nn.
\eea

So we have proven that $\Gamma$ is the heat kernel on $\bR^2\times[-\pi,\pi]$ for those functions whose integrals on $S^1$ are zero. Now we define
\[
G(x,y)=\int^\i_0\G(t;x,y)dt,
\]which is finite except when $x=y$, due the exponential decaying property of $\Gamma_1$. See the next lemma for details.
Let $f$ be a smooth, compactly supported function on $\bR^2 \times S^1$, whose integral on $S^1$ is $0$. From the above computation, it is easy to see that
\bea
&&-\int \Dl_x G(x,y)f(y) dy= \int \int^\i_0-\Dl_x\G(t; x, y) dt f(y) dy \nn\\
&&\qq\qq=-\int \int^\infty_0 \p_t\G(t; x, y) dt f(y) dy =\int \G(t;x,y)|_{t=0} f(y) dy=f(x).\nn
\eea Hence
\[
\Delta [u(x) - \int  G(x,y)f(y) dy]=0.
\]Note the function $u - \int  G(x,y)f(y) dy$ is bounded. Therefore, the classical Liouville theorem implies that
\[
u= \int  G(x,y)f(y) dy + C.
\]
 This shows that $G$ is the Green function on $\bR^2\times[-\pi,\pi]$ of those functions with its integral on $S^1$ is zero. \qed

\begin{lemma}\label{l10.1}
Let $G(x,y)$ be the Green function on $\bR^2\times S^1$ defined above. Then we have the following estimates: for a positive constant $c_0$,
\be
|G(x,y)|\ls \f{1}{|x-y|}e^{-c_0|x'-y'|},\q |\na G(x,y)|\ls \f{1}{|x-y|^2}e^{-c_0|x'-y'|}.
\ee
with $x'=(x_1,x_2)$ and $y'=(y_1,y_2)$.
\end{lemma}

\pf   We will prove Lemma \ref{l10.1}  by dividing the estimate in two cases: $|x'-y'|> 1$ and $|x'-y'|\leq 1$.

\noindent\textbf{Case1: $\boldsymbol{|x'-y'|> 1}$}, we will show
\be\label{4.6}
{|G(x,y)|+|\na G(x,y)|\ls e^{-\f{|x'-y'|}{2}}}.
\ee
From \eqref{4.1e}, we have
\be\label{10.5}
\begin{aligned}
G(x,y)&=\int^\i_0(4\pi t)^{-1}e^{-\f{|x'-y'|^2}{4t}}\f{1}{\pi}\sum_{m\in\bZ^+} e^{-m^2 t} \cos (m(x_3-y_3))dt\\
&=\f{1}{4\pi^2}\sum^\i_{m=1}\int^\i_0 t^{-1}e^{-\f{|x'-y'|^2}{4t}}e^{-m^2t}dt\cos (m(x_3-y_3))\\
&:=\f{1}{4\pi^2}\sum^\i_{m=1} I_m\cos (m(x_3-y_3)).
\end{aligned}
\ee

Now we come to estimate $I_m$. Making a variable change, it is easy to see that
\be\label{10.6}
\begin{aligned}
I_m&=\int^\i_0 t^{-1}e^{-\f{(m|x'-y'|)^2}{4t}}e^{-t}dt\\
&=\Big(\int^{\f{m|x'-y'|}{2}}_0+\int^\i_{\f{m|x'-y'|}{2}}\Big) t^{-1}e^{-\f{(m|x'-y'|)^2}{4t}}e^{-t}dt\\
&=2\int^\i_{\f{m|x'-y'|}{2}} t^{-1}e^{-\f{(m|x'-y'|)^2}{4t}}e^{-t}dt\\
&\ls \int^\i_{\f{m|x'-y'|}{2}} t^{-1}e^{-t}dt\\
&\ls \ln \big(\f{1}{m|x'-y'|}+2\big)e^{-\f{m|x'-y'|}{2}}.
\end{aligned}
\ee
From \eqref{10.5} and \eqref{10.6}, we have
\be\label{10.7}
\begin{aligned}
\big| G(x,y)\big |&\ls\sum^\i_{m=1}\ln \big(\f{1}{m|x'-y'|}+2\big)e^{-\f{m|x'-y'|}{2}}\\
&\ls \ln \big(\f{1}{|x'-y'|}+2\big)\sum\limits_me^{-\f{m|x'-y'|}{2}}.\\
&\ls \ln \big(\f{1}{|x'-y'|}+2\big)\f{e^{-\f{|x'-y'|}{2}}}{1-e^{-\f{|x'-y'|}{2}}}\\\
&\ls e^{-\f{|x'-y'|}{2}}.
\end{aligned}
\ee

This finishes the proof for $G(x,y)$.

Also from \eqref{10.5}, we have
\be\label{10.8}
\begin{aligned}
\big|\p_{x_3} G(x,y)\big|&=\Big|\f{1}{4\pi^2}\sum^\i_{m=1}m\int^\i_0 t^{-1}e^{-\f{|x'-y'|^2}{4t}}e^{-m^2t}dt\cos' (m(x_3-y_3))\Big|\\
&\ls \Big|\sum^\i_{m=1}\int^\i_0 t^{-\f{3}{2}}e^{-\f{|x'-y'|^2}{4t}}e^{-\f{m^2t}{2}}dt\Big|\\
&\ls \Big|\sum^\i_{m=1}m\int^\i_0 t^{-\f{3}{2}}e^{-\f{(m|x'-y'|)^2}{4t}}e^{-\f{t}{2}}dt\Big|\\
&\ls \Big|\sum^\i_{m=1}m\Big(\int^\i_{\f{m|x'-y'|}{2}}+\int^{\f{m|x'-y'|}{2}}_0\Big) t^{-\f{3}{2}}e^{-\f{(m|x'-y'|)^2}{4t}}e^{-\f{t}{2}}dt\Big|\\
&\ls \f{1}{|x'-y'|}\Big|\sum^\i_{m=1}\int^\i_{\f{m|x'-y'|}{2}}(t^{-1}+ t^{-\f{1}{2}})e^{-\f{t}{2}}dt\Big|\\
&\ls e^{-\f{|x'-y'|}{2}}.
\end{aligned}
\ee

The estimate of $\p_{y_3} G(x,y)$ is the same as \eqref{10.8}.

From \eqref{10.5}, we have
\be\label{10.11}
\begin{aligned}
\big|\p_{x',y'} G(x,y)\big|&\ls \Big|\sum^\i_{m=1}\int^\i_0 t^{-1}\f{|x'-y'|}{t}e^{-\f{|x'-y'|^2}{4t}}e^{-{m^2t}}dt\Big|\\
&\ls |x'-y'|\Big|\sum^\i_{m=1}m^2\int^\i_0 t^{-2}e^{-\f{(m|x'-y'|)^2}{4t}}e^{-{t}}dt\Big|\\
&\ls |x'-y'|\Big|\sum^\i_{m=1}m^2\Big(\int^\i_{\f{m|x'-y'|}{2}}+\int^{\f{m|x'-y'|}{2}}_0\Big) t^{-2}e^{-\f{(m|x'-y'|)^2}{4t}}e^{-{t}}dt\Big|\\
&\ls\f{1}{|x'-y'|}\Big|\sum^\i_{m=1}\int^\i_{\f{m|x'-y'|}{2}}(t^{-1}+1)e^{-t}dt\Big|\\
&\ls e^{-\f{|x'-y'|}{2}}.
\end{aligned}
\ee

\noindent\textbf{Case2: $\boldsymbol{|x'-y'|\leq 1}$}, we will show there exists a positive constant $c$ such that
\be\label{4.12}
{|G(x,y)|\leq \f{1}{|x-y|}e^{-c|x-y|},\q |\na G(x,y)|\leq \f{1}{|x-y|^2}e^{-c|x-y|}}.
\ee

Let $\G(t,x,y)$ be the heat kernel in \eqref{ee4.2} on $\bR^2\times S^1$. The explicit formula for the heat kernel in $\bR^2$ together with
 the well known global behavior of the heat kernel on $S^1$ and gradient with mean $0$  tell us that there exist two positive constants $c_1$ and $c_2$ such that
\be\label{4.13}
|\G(t, x,y)|\ls t^{-3/2}e^{-c_1\f{|x-y|^2}{t}}e^{-c_2t},\q |\na \G(x,y,t)|\ls t^{-2}e^{-c_1\f{|x-y|^2}{t}}e^{-c_2t}.
\ee In fact one can choose $c_2$ as any positive number less than $1/4$.
Readers can see, for example \cite{TZ:1} for more details for the estimate of \eqref{4.13}, even when $S^1$ is replaced by more general compact manifold. Integrating \eqref{4.13} on $t$ from $0$ to $+\i$ implies \eqref{4.12}.\\

 By choosing $c_0$ small ($c_0=\f{1}{3}\min\{c,\f{1}{2}\}$) and remembering the fact $|x_3-y_3|\leq 4\pi$, we see that \eqref{4.5} is a direct consequence of
\eqref{4.6} and \eqref{4.12}.\ef

\begin{lemma}\label{l10.3}
Let ${G}(x,y)$ be defined as above. Denote $x=(r\cos\th,r\sin\th,z)$ and $y=(\rho\cos\phi,\rho\sin\phi,l)$, then we have the following estimates.

For $|\rho-r|\leq \f{1}{4}r$,
\be\label{10.15}
\int^{2\pi}_0|\na {G}(x,y)|d\phi\ls \f{1}{\rho(|\rho-r|+|z-l|)}e^{-c_0|\rho-r|}.
\ee
\be\label{e4.13e}
\int^{2\pi}_0|{G}(x,y)|d\phi\ls e^{-c_0|\rho-r|}\f{1}{r}\ln \Big(2+\f{r}{|\rho-r|}\Big).
\ee
\end{lemma}
When $\f{1}{8}r \le |\rho-r|\le \f{1}{4}r$, we have
\be\label{4.15}
\int^{2\pi}_0\Big(|{G}(x,y)|+|\na {G}(x,y)|\Big)d\phi\ls e^{-c_0|\rho-r|}.
\ee
\pf
Remember that $x'=(r\cos\th,r\sin\th), y'=(\rho\cos\phi,\sin\phi)$, then we can get
\be
|x'-y'|=\s{\rho^2+r^2-2\rho r\cos(\th-\phi)}\xlongequal{\th=0}\s{(\rho-r)^2+4\rho r\sin^2\f{\phi}{2}}.\nn
\ee

Denote $J:=\int^{2\pi}_0\f{1}{|x-y|}d\phi$ and $K:=\int^{2\pi}_0\f{1}{|x-y|^2}d\phi$. Next we will prove
\be\label{10.18}
J\ls \f{1}{r}\ln (2+\f{r}{|\rho-r|}),\q K\ls \f{1}{\rho(|\rho-r|+|z-l|)}.
\ee
Then \eqref{4.5} and \eqref{10.18} together prove \eqref{10.15}, \eqref{e4.13e}.

It is easy to see that
\be\label{10.19}
\begin{aligned}
J&\ls\int^{2\pi}_0\f{1}{|x'-y'|}d\phi\\
&\ls\int^{\pi}_0 \f{d\phi}{\s{(\rho-r)^2+4\rho r\sin^2\phi}}\\
&\ls\int^{\f{\pi}{2}}_0 \f{d\phi}{\s{(\rho-r)^2+4\rho r\sin^2\phi}}\\
& \thickapprox \f{1}{r}\int^{\f{\pi}{2}}_0 \f{d\phi}{\s{\kappa^2+\sin^2\phi}} \q \text{with}\q \kappa^2=\f{|\rho-r|^2}{4\rho r}\ls 1\\
& \thickapprox \f{1}{r}\Big(\int^{\f{\pi}{4}}_0+\int^{\f{\pi}{2}}_{\f{\pi}{4}}\Big) \f{d\phi}{\s{\kappa^2+\sin^2\phi}}.
\end{aligned}
\ee
If $\phi\in [0,\f{\pi}{4}]$, $\sin\phi\thickapprox \phi$ and if $\phi\in [\f{\pi}{4},\f{\pi}{2}]$, $\sin\phi\thickapprox 1$, then \eqref{10.19} implies that
\bea
&&J\thickapprox \f{1}{r}\Big (\int^{\f{\pi}{4}}_0\f{d\phi}{\s{\kappa^2+\phi^2}}+\int^{\f{\pi}{2}}_{\f{\pi}{4}}\f{d\phi}{\s{\kappa^2+1}}\Big)\nn\\
&&\q \thickapprox \f{1}{r}\Big (\int^{\f{\pi}{4\kappa}}_0\f{d\phi}{1+\phi}+1\Big)\nn\\
&&\q\ls \f{1}{r}\ln\Big(2+\f{1}{\kappa}\Big)\nn\\
&&\q\ls \f{1}{r}\ln\Big(2+\f{r}{|\rho-r|}\Big)\nn\\
&&\q\ls \f{1}{r}\ln\Big(2+\f{r}{|\rho-r|}\Big).\nn
\eea
The estimate for $K$ in \eqref{10.18} will essentially the same as $J$, we omit the detail. Then \eqref{4.15} is a direct consequence of \eqref{10.15}, \eqref{e4.13e} when $\f{1}{8}r\leq|\rho-r|\leq \f{1}{4}r$. \ef

\subsection{First decay on $u$ and $ w$ with the assumption $\int^\pi_{-\pi}u^\th dz=\int^\pi_{-\pi}u^z dz=0$}

We start with a simple

{\it \noindent Observation.

Assume $u=u^re_r+u^\th e_\th+u^ze_z$ is a z-periodic axisymmetric solution of \eqref{e10.1}, then for $r\neq 0$ and $k\in\bN$, we have
\be\label{10.10}
\int^\pi_{-\pi} \p^k_ru^r(r,z)dz=0.
\ee
}

This can be seen from the formula $u^r = - \p_z L^\theta $ where $L^\theta$ is the angular stream function which is also periodic in $z$ variable.

Alternatively, by the incompressible condition, we have
\be
\f{1}{r}\p_r(ru^r)+\p_zu^z=0.\nn
\ee
Integrating the equation on $z$ from $-\pi$ to $\pi$ implies
\be
\f{1}{r}\p_r\big(r\int^\pi_{-\pi} u^rdz\big)=-\int^\pi_{-\pi} \p_z u^zdz=u^z\Big|^\pi_{-\pi}=0\nn
\ee
which indicates
\be
\p_r\big(r\int^\pi_{-\pi} u^rdz\big)=0.\nn
\ee

Thus we have
\be
\int^\pi_{-\pi} u^rdz=\f{1}{r}\Big(\int^r_0\p_{\t{r}}\big(\t{r}\int^\pi_{-\pi} u^rdz\big)d\t{r}\Big)=0.\nn
\ee
Differentiating the above equation $k$ times with respect to $r$, we prove \eqref{10.10}.
The rest of the subsection is divided into 3 steps.

\medskip

{\bf Step 1.} First decay of $w^\th$.

Again we pick a point $x_0 \in \bR^3$ such that $|x'_0| \equiv \la $ is large and carry out the
scaling for the velocity and vorticity:
\bea
&&\t{u}(\t{x})=\la u(\la \t{x})=\la u(x)\nn\\
&&\t{w}(\t{x})=\la^2w(\la \t{x})=\la^2w(x)\nn
\eea
where $\t{x}=\f{x}{\la}$.

For simplification of notation, we temporarily drop the $``\sim"$ symbol when computations take place under the scaled sense. Define the domains
\be
\mathcal{D}_1=\{(r,\th,z):\f{1}{2}<r<\f{3}{2},\ 0\leq \th\leq 2\pi,\ z\in[-\f{\pi}{\la},\f{\pi}{\la}]\}\nn
\ee
and
\be
\mathcal{D}_2=\{(r,\th,z):\f{3}{4}<r<\f{5}{4},\ 0\leq \th\leq 2\pi,\ z\in[-\f{\pi}{\la},\f{\pi}{\la}]\}.\nn
\ee
Almost identical computations to \eqref{4.4} and \eqref{4.5} lead to
\be\label{e4.20}
\begin{aligned}
&\q\ \|\big(\nabla w^r,\nabla w^z\big)\|^2_{L^2(\mathcal{D}_2)}\\
&\leq C(1+\|(u^r,u^z)\|^2_{L^\i(\mathcal{D}_1)})\|(w^r,w^z)\|^2_{L^2(\mathcal{D}_1)},
\end{aligned}
\ee
and
\be\label{e4.21}
\begin{aligned}
&\q\ \|\nabla w^\th\|^2_{L^2(\mathcal{D}_2)}\\
&\leq C \big(1+\|(u^r,u^\th,u^z)\|_{L^\i(\mathcal{D}_1)}\big)\|(w^r,w^\th)\|^2_{L^2(\mathcal{D}_1)}.
\end{aligned}
\ee Comparing with the full space case in the previous section,   since all functions  here are periodic in $z$ variable, the cut off function that leads to the above inequality  only depends on $r$.

Recall the Brezis-Gallouet inequality
of Lemma \ref{l3.1}.
Let $f\in H^2(\O)$ where $\O \subset \bR^2$ is a piecewise smooth domain. Then there exists a constant $C_{\O}$, depending only on $\O$, such that
\be
\|f\|_{L^\i(\O)}\leq C_{\O}\|f\|_{H^1(\O)}\ln^{1/2} \left(e+\f{\|\Dl f\|_{L^2(\O)}}{\|f\|_{H^1(\O)}} \right).\nn
\ee
Now we select the 2 dimensional domain
\be
\bar{\mathcal{D}}_2:=\{(r,z):\f{3}{4}<r<\f{5}{4},\ z\in[-\f{\pi}{\la},\f{\pi}{\la}]\}.\nn
\ee
If we choose $\O=\bar{\mathcal{D}}_2$ and by scaling in the $z$ direction only, we  have
\be
C_{\bar{\mathcal{D}}_2}\leq C\la^{1/2}
\ee where $C$ is independent of $\bar{\mathcal{D}}_2$ and $\la$. i.e.
 If $f\in H^2(\bar{\mathcal{D}}_2)$, then there exists an absolute constant $C$, such that
\be
\label{BGcbar}
\|f\|_{L^\i(\bar{\mathcal{D}}_2)}\leq C \la^{1/2} \|f\|_{H^1(\bar{\mathcal{D}}_2)}\ln^{1/2} \big(e+\f{\|\Dl f\|_{L^2(\bar{\mathcal{D}}_2)}}{\|f\|_{H^1(\bar{\mathcal{D}}_2)}}\big).
\ee Notice that the inequality is worse than the corresponding one in the full space case by
a factor of $\la^{1/2}$.

Using \eqref{e4.21} and \eqref{BGcbar} , we deduce
\be\label{e4.22}
\begin{aligned}
&\q\ \|w^\th\|_{L^\i(\bar{\mathcal{D}}_2)}  \\
&\ls \la^{1/2}\|w^\th\|_{H^1(\bar{\mathcal{D}}_2)}\ln^{1/2} \big(e+\f{\|\Dl w^\th\|_{L^2(\bar{\mathcal{D}}_2)}}{\|w^\th\|_{H^1(\bar{\mathcal{D}}_2)}}\big)\\
&\ls \la^{1/2}\big(1+\|w^\th\|_{H^1(\bar{\mathcal{D}}_2)}\big)\ln^{1/2} \big(e+\|\Dl w^\th\|_{L^2(\bar{\mathcal{D}}_2)}\big)\\
&\ls  \la^{1/2}\Big(1+(1+\|(u^r,u^\th,u^z)\|^{1/2}_{L^\i(\mathcal{D}_1)})\|w^r,w^\th\|_{L^2({\mathcal{D}}_1)}\Big)\\
&\q\ \times\ln^{1/2} \big(e+\|\Dl w^\th\|_{L^2({\mathcal{D}}_1)}\big);
\end{aligned}
\ee
and  using \eqref{e4.20} and \eqref{BGcbar}, we find:
\be\label{e4.23}
\begin{aligned}
&\q\ \|(w^r,w^z)\|_{L^\i(\bar{\mathcal{D}}_2)}  \\
&\ls \la^{1/2}\|(w^r,w^z)\|_{H^1(\bar{\mathcal{D}}_2)}\ln^{1/2} \big(e+\f{\|(\Dl w^r,\Dl w^z)\|_{L^2(\bar{\mathcal{D}}_2)}}{\|(w^r, w^z)\|_{H^1(\bar{\mathcal{D}}_2)}}\big)\\
&\ls \la^{1/2}\big(1+\|(w^r,w^z)\|_{H^1(\bar{\mathcal{D}}_2)}\big)\ln^{1/2} \big(e+\|(\Dl w^r,\Dl w^z)\|_{L^2(\bar{\mathcal{D}}_2)}\big)\\
&\ls \la^{1/2}\Big(1+(1+\|(u^r,u^z)\|_{L^\i(\mathcal{D}_1)})\|(w^r,w^z)\|_{L^2({\mathcal{D}}_1)}\Big)\\
&\q\ \times\ln^{1/2} \big(e+\|(\Dl w^r,\Dl w^z)\|_{L^2({\mathcal{D}}_1)}\big).
\end{aligned}
\ee

Now scaling back to the domains
\be
\mathcal{D}_{1, \la}=\{(r,\th,z):\f{\la}{2}<r<\f{3 \la}{2},\ 0\leq \th\leq 2\pi,\ z\in[-\pi, \pi]\},\nn
\ee
and
\be
\mathcal{D}_{2, \la}=\{(r,\th,z):\f{3\la}{4}<r<\f{5\la}{4},\ 0\leq \th\leq 2\pi,\ z\in[-\pi, \pi]\},\nn
\ee
 we deduce

\be\label{e4.24}
\begin{aligned}
&\q\ \la^2\|w^\th\|_{L^\i(\mathcal{D}_{2,\la})}\\
&\ls\la^{1/2}\Big(1+\big(1+\la^{1/2}\|(u^r,u^\th,u^z)\|^{1/2}_{L^\i(\mathcal{D}_{1,\la})}\big)\la^{1/2}\|(w^r,w^\th)\|_{L^2(\mathcal{D}_{1,\la})}\Big)\\
&\q\ \times\ln^{1/2} \big(\la^{5/2}\|\Dl w^\th\|_{L^2(\mathcal{D}_{2,\la})}+e\big)\\
&\ls\la\Big(1+\big(1+\la^{1/2}\|(u^r,u^\th,u^z)\|^{1/2}_{L^\i(\mathcal{D}_{1,\la})}\big)\|(w^r,w^\th)\|_{L^2(\mathcal{D}_{1,\la})}\Big)\\
&\q\ \times\ln^{1/2} \big(\la^{5/2}\|\Dl w^\th\|_{L^2(\mathcal{D}_{2,\la})}+e\big)
\end{aligned}
\ee
and
\be\label{e4.25}
\begin{aligned}
&\q\ \la^2\|(w^r,w^z)\|_{L^\i(\mathcal{D}_{2,\la})}\\
&\leq C\la\Big(1+\big(1+\la\|(u^r,u^z)\|_{L^\i(\mathcal{D}_{1,\la})}\big)\|(w^r,w^z)\|_{L^2(\mathcal{D}_{1,\la})}\Big)\\
&\q\ \times\ln^{1/2} \big(\la^{5/2}\|(\Dl w^r,\Dl w^z)\|_{L^2(\mathcal{D}_{2,\la})}+e\big).
\end{aligned}
\ee

Also scaling back of \eqref{e4.20} and \eqref{e4.21} indicate that

\be\label{e4.26}
\begin{aligned}
&\q\ \|\big(\nabla w^r,\nabla w^z\big)\|_{L^2(\mathcal{D}_{2,\la})}\\
&\ls\la^{-1}(1+\la\|(u^r,u^z)\|_{L^\i(\mathcal{D}_{1,\la})})\|(w^r,w^z)\|_{L^2(\mathcal{D}_{1,\la})},
\end{aligned}
\ee
and
\be\label{e4.27}
\begin{aligned}
&\q\ \|\nabla w^\th\|_{L^2(\mathcal{D}_{2,\la})}\\
&\ls \la^{-1} \big(1+\la^{1/2}\|(u^r,u^\th,u^z)\|^{1/2}_{L^\i(\mathcal{D}_{1,\la})}\big)\|(w^r,w^\th)\|_{L^2(\mathcal{D}_{1,\la})}.
\end{aligned}
\ee

From \eqref{e4.24} and assumption (\ref{dr2r}), we reach the first decay estimate of $w^\th$
\be\label{e4.28}
|w^\th(x)|\leq C_0\big(\f{\ln r}{r}\big)^{1/2}.
\ee
\medskip

{\bf Step 2.} First decay of $u$ and $(w^r, w^z)$.

Using the Biot-Savart law, for a cutoff functions $\psi=\psi(x')$, which is independent of $x_3$, we know, for any smooth, divergence free vector field $b$, that
\be
-\Dl(b \psi)=\psi\na\times (\na \times b) -2\na\psi\cdot\na b -b\Dl\psi.\nn
\ee

Then using Green function on $\bR^2\times S^1$, we have
\be\label{10.13}
\begin{aligned}
b(x)&=\int_{S^1}\int_{\bR^2}G(x,y)\psi\na\times (\na \times b) dy\\
&-2\int_{S^1}\int_{\bR^2}G(x,y)\na\psi\cdot\na b dy-\int_{S^1}\int_{\bR^2}G(x,y)(\Dl\psi) b dy.
\end{aligned}
\ee

If we write $b=u^re_r+u^ze_z$, then $\na\times b=w^\th e_\th$.  Let $x=(r\cos\th,r\sin\th,z)$, $y=(\rho\cos\phi,\rho\sin\phi,l)$. Then from \eqref{10.13}, we have
\bea
&&u^r(x)=\int_{S^1}\int_{\bR^2}G(x,y)\psi(\na\times (w^\th e_\th))\cdot e_rdy\nn\\
&&\qq-2\int_{S^1}\int_{\bR^2}G(x,y)(\na\psi\cdot\na b)\cdot e_rdy-\int_{S^1}\int_{\bR^2}G(x,y)(\Dl\psi) b\cdot e_r d y\nn\\
&&\q\ =\int^{\pi}_{-\pi}\int^{2\pi}_0\int^\i_0G(x,y)\psi \p_l w^\phi \cos(\phi-\th)\rho d\rho d\phi dl\nn\\
&&\qq -2\int^{\pi}_{-\pi}\int^{2\pi}_0\int^\i_0G(x,y)\p_\rho\psi\p_\rho u^\rho \cos(\phi-\th) \rho d\rho d\phi dl\nn\\
&&\qq -\int^{\pi}_{-\pi}\int^{2\pi}_0\int^\i_0G(x,y)(\p^2_\rho\psi+\f{1}{\rho}\p_\rho\psi) u^\rho \cos(\phi-\th) \rho d\rho d\phi dl.\nn
\eea

Since our solution is axisymmetric, we can set $\th=0$. Also due to the fact that $\psi$ is independent of $l$ and the fact that
\be
\int^\pi_{-\pi}\p_l w^\th dl=\int^\pi_{-\pi}\p_\rho u^\rho dl=\int^\pi_{-\pi} u^\rho dl=0, \nn
\ee we get

\be\label{10.14}
\begin{aligned}
u^r&=\int^{\pi}_{-\pi}\int^{2\pi}_0\int^\i_0 {G}(x,y)\psi \p_l w^\phi \cos \phi\rho d\rho d\phi dl\\
&\q\ -2\int^{\pi}_{-\pi}\int^{2\pi}_0\int^\i_0 {G}(x,y)\p_\rho\psi\p_\rho u^\rho \cos\phi \rho d\rho d\phi dl\\
&\q\ -\int^{\pi}_{-\pi}\int^{2\pi}_0\int^\i_0 {G}(x,y)(\p^2_\rho\psi+\f{1}{\rho}\p_\rho\psi) u^\rho \cos\phi \rho d\rho d\phi dl.\\
&=\underbrace{-\int^{\pi}_{-\pi}\int^\i_0\Big(\int^{2\pi}_0\p_l {G}(x,y)\cos \phi d\phi\Big)\psi  w^\phi \rho d\rho  dl}_{I_1}\\
&\q\ -\underbrace{2\int^{\pi}_{-\pi}\int^\i_0\Big(\int^{2\pi}_0 {G}(x,y)\cos\phi d\phi\Big)\p_\rho\psi\p_\rho u^\rho  \rho d\rho dl}_{I_2}\\
&\q\ -\underbrace{\int^{\pi}_{-\pi}\int^\i_0\Big(\int^{2\pi}_0 {G}(x,y)\cos\phi d\phi\Big)(\p^2_\rho\psi+\f{1}{\rho}\p_\rho\psi) u^\rho \rho d\rho d\phi dl}_{I_3}.
\end{aligned}
\ee

Now we need to estimate $I_1,I_2,I_3$ separately. From Lemma \ref{l10.3},
\bea
&&|I_1|\ls \int^\pi_{-\pi}\int_{|\rho-r|\leq \f{1}{4}r}\f{1}{\rho}\f{1}{|\rho-r|+|z-l|}e^{-c_0|\rho-r|}|\psi||w^\phi|\rho d\rho dl\nn\\
&&\q \ls \sup\limits_{(\rho,l)\in[\f{3}{4}r,\f{5}{4}r]\times[-\pi,\pi]} |w^\phi(\rho,l)|\int^\pi_{-\pi}\int_{|\rho-r|\leq \f{1}{4}r} \f{1}{|\rho-r|+|z-l|}e^{-c_0|\rho-r|}d\rho dl\nn\\
&&\q \ls \big(\f{\ln r}{r}\big)^{1/2}\underbrace{\int_{|\rho-r|\leq \f{1}{4}r}\ln\big(1+\f{\pi}{|\rho-r|}\big)e^{-c_0|\rho-r|} d\rho.}_{J}\nn
\eea
When $r$ is large, it is easy to see that $J\ls 1$. So, we can get
\be\label{10.20}
|I_1|\ls r^{-1/2}(\ln r)^{1/2}.
\ee

For $I_2$, using Lemma \ref{l10.3}, we have
\be\label{10.21}
\begin{aligned}
|I_2|&\ls \int^\pi_{-\pi}\int_{\f{1}{8}r\leq|\rho-r|\leq \f{1}{4}r }e^{-c_0|\rho-r|}|\p_\rho \psi||\p_\rho u^\rho|\rho d\rho dl\\
&\ls \sup |\p_\rho u^\rho|\int_{\f{1}{8}r\leq|\rho-r|\leq \f{1}{4}r }e^{-c_0|\rho-r|} d\rho \\
&\ls re^{-\f{c_0}{8}r}\ls e^{-\f{c_0}{16}r}.
\end{aligned}
\ee
We can also prove that
\be\label{10.22}
|I_3|\ls e^{-\f{c_0}{16}r}.
\ee
At last, from \eqref{10.20},\eqref{10.21} and \eqref{10.22}, we can get for large $r$
\be
\lab{ur1decay}
|u^r|\ls r^{-1/2}(\ln r)^{1/2}.
\ee

Using \eqref{10.13} in cylindrical coordinates, setting $x=(r\cos\th,\r\sin\th,z),y=(\rho\cos\phi,\rho\sin\phi,l)$ and using integration by parts, we can get
\bea
&&u^z=\int^{\pi}_{-\pi}\int^{2\pi}_0\int^\i_0G(x,y)\psi\p_\rho(\rho w^\phi) d\rho d\phi dl\nn\\
&&\qq -2\int^{\pi}_{-\pi}\int^{2\pi}_0\int^\i_0G(x,y)\p_\rho\psi\p_\rho u^z \rho d\rho d\phi dl\nn\\
&&\qq -\int^{\pi}_{-\pi}\int^{2\pi}_0\int^\i_0G(x,y)(\p^2_\rho\psi+\f{1}{\rho}\p_\rho\psi) u^z \rho d\rho d\phi dl\nn
\eea

Since $w^\th=\p_z u^r-\p_r u^z$  and our assumption that the integral of $u^z$ with respect to $z$ on $[-\pi,\pi]$, we can get
\be\label{10.24}
\int^\pi_{-\pi}w^\th dz=0.
\ee
Due to \eqref{10.24} and axisymmetry of $u$, we can set $\th=0$. Then using integration by parts, we obtain
\be\label{e4.36}
\begin{aligned}
u^z&=\int^{\pi}_{-\pi}\int^\i_0\Big(\int^{2\pi}_0 {G}(x,y) d\phi\Big)\psi\p_\rho(\rho w^\phi) d\rho dl\\
&\q\ -2\int^{\pi}_{-\pi}\int^\i_0\Big(\int^{2\pi}_0 {G}(x,y) d\phi\Big)\p_\rho\psi\p_\rho u^z \rho d\rho dl\\
&\q\ -\int^{\pi}_{-\pi}\int^\i_0\Big(\int^{2\pi}_0{G}(x,y) d\phi\Big)(\p^2_\rho\psi+\f{1}{\rho}\p_\rho\psi) u^z \rho d\rho dl.\\
&=\underbrace{-\int^{\pi}_{-\pi}\int^\i_0\Big(\int^{2\pi}_0\p_\rho {G}(x,y) d\phi\Big)\psi w^\phi \rho d\rho dl}_{J_1}\\
&\q\ \underbrace{-\int^{\pi}_{-\pi}\int^\i_0\Big(\int^{2\pi}_0 {G}(x,y) d\phi\Big)\p_\rho\psi w^\phi \rho d\rho dl}_{J_2}\\
&\q\ \underbrace{-2\int^{\pi}_{-\pi}\int^\i_0\Big(\int^{2\pi}_0 {G}(x,y) d\phi\Big)\p_\rho\psi\p_\rho u^z \rho d\rho dl}_{J_3}\\
&\q\ \underbrace{-\int^{\pi}_{-\pi}\int^\i_0\Big(\int^{2\pi}_0{G}(x,y) d\phi\Big)(\p^2_\rho\psi+\f{1}{\rho}\p_\rho\psi) u^z \rho d\rho dl}_{J_4}.
\end{aligned}
\ee

Then the estimates of $J_1$ will essentially be the same as $I_1$ and estimates of $J_2, J_3, J_4$, will essentially be the same as $I_2,I_3, I_4$. At last, we can get
\be\label{10.27}
|u^z|\ls \f{1}{r^{1/2}}(\ln r)^{1/2}.
\ee

Using the decay of $u^r$ and $u^z$ so far, i.e. substituting \eqref{ur1decay} and \eqref{10.27} into the right hand side of \eqref{e4.25} gives the first decay of $(w^r, w^z)$:
\be
\lab{wrwz1decay}
|(w^r, w^z)|\ls \f{\ln r}{r^{1/2}}.
\ee

Since $\int^\pi_{-\pi}u^\th dz=0$ by assumption, there exists, for any $r$, a $z_0$ such that $u^\th(r,z_0)=0$ where $z_0$ may be different for different $r$.  Then, using $\p_z u^\th= - w^r$, we find
\be\label{e4.39}
|u^\th(r,z)|=|\int^z_{z_0}\p_z u^\th dz|\ls |\p_z u^\th|_{L^\i}2\pi\ls \f{1}{r^{1/2}}\ln r.
\ee
\medskip
{\bf Step 3.} Almost first order decay of $u$ and $w$ by iteration.

Now we can iterate and improve the decay of $u$ by using  \eqref{e4.24},\eqref{e4.25},\eqref{10.14},\eqref{e4.36} and \eqref{e4.39} with the following order.

applying first decay of $w^\th$ on  \eqref{10.14},\eqref{e4.36} $\Rightarrow$  decay of $u^r,u^z$;

applying decay of $u^r, u^z$ on \eqref{e4.25} $\Rightarrow$ decay of $ w^r, w^z$;

applying decay of $ w^r, w^z$ on \eqref{e4.39} $\Rightarrow$ decay of $u^\th$;

applying decay of $u^r, u^z, u^\th$ on (\ref{e4.24}) $\Rightarrow$  second decay of $w^\th$, and so on.

\noindent The detail is as follows. Define the set
\be
S_n:=\big\{(\rho,z)\big | \{|\rho-r|\leq \f{r}{2^{n+2}}\}\times[-\pi,\pi]\big\}.\nn
\ee

From \eqref{10.14} and \eqref{e4.36}, we have
\be
\|(u^r,u^z)\|_{L^\i(S_{n+1})}\leq C_n\|w^\th\|_{L^\i(S_n)}+C_ne^{-\f{c_0}{16}r},\nn
\ee
From \eqref{e4.25}, we have
\be
\|(w^r,w^z)\|_{L^\i(S_{n+2})}\leq C_{n+1}\|(u^r,u^z)\|_{L^\i(S_{n+1})}(\ln r)^{1/2},\nn
\ee
From \eqref{e4.39}, we have
\be
\|u^\th\|_{L^\i(S_{n+2})}\leq 2\pi\|w^r\|_{L^\i(S_{n+2})}.\nn
\ee
From \eqref{e4.24}, we have
\be
\|w^\th\|_{L^\i(S_{n+3})}\leq C_{n+2}\|(u^r,u^\th,u^z)\|^{1/2}_{L^\i(S_{n+2})}r^{-1/2}(\ln r)^{1/2}.\nn
\ee
where $\lim\limits_{n\rightarrow +\i}C_{n}=+\i $.

 Recall from step 1 that $\|w^\th\|_{L^\i(S_{0})}\leq{C}_0r^{-1/2}(\ln r)^{1/2}$. After iteration $n$ times,  we can get
\be\label{e4.39e}
\begin{aligned}
\|w^\th\|_{L^\i(S_{3n})}&\leq A_{n}r^{-1+\f{1}{2^{n+1}}}(\ln r)^{\f{3}{2}-\f{1}{2^{n-1}}}\\
&\q\ +A_nr^{-1/2}(\ln r)^{3/4} e^{-\f{c_0}{32}r}\\
&\leq 2A_nr^{-1+\f{1}{2^{n+1}}}(\ln r)^{\f{3}{2}-\f{1}{2^{n-1}}}.
\end{aligned}
\ee
And also we have
\be\label{e4.40e}
\begin{aligned}
\|(u^r,u^z)\|_{L^\i(S_{3n+1})}&\leq C_{3n}\|w^\th\|_{L^\i(S_{3n})}+C_{3n}e^{-\f{c_0}{16}r},\\
&\ls B_nr^{-1+\f{1}{2^{n+1}}}(\ln r)^{\f{3}{2}-\f{1}{2^{n-1}}},
\end{aligned}
\ee
\be\label{e4.41e}
\begin{aligned}
\|(w^r,w^z)\|_{L^\i(S_{3n+2})}&\leq C_{3n+1}\|(u^r,u^z)\|_{L^\i(S_{3n+1})}(\ln r)^{1/2},\\
&\leq D_nr^{-1+\f{1}{2^{n+1}}}(\ln r)^{2-\f{1}{2^{n-1}}},
\end{aligned}
\ee
\be\label{e4.42e0}
\begin{aligned}
\|u^\th\|_{L^\i(S_{3n+2})}&\leq 2\pi\|w^r\|_{L^\i(S_{3n+2})},\\
&\leq 2\pi D_nr^{-1+\f{1}{2^{n+1}}}(\ln r)^{2-\f{1}{2^{n-1}}},
\end{aligned}
\ee
where the constants $A_n, B_n, D_n$ satisfy
\be
\lim\limits_{n\rightarrow +\i}A_n=\lim\limits_{n\rightarrow +\i}B_n=\lim\limits_{n\rightarrow +\i}D_n=+\i.\nn
\ee

Consider the domain
\be
\Omega_\dl:=\big\{(\rho,z)\big | \{|\rho-r|\leq \f{\dl^3}{4}r\}\times[-\pi,\pi]\big\}.\nn
\ee
Then for large $r$, we can get
\be\label{e4.42e}
\|(u^r,u^\th,u^z)\|_{L^\i(\Omega_\dl)}\leq C_\dl r^{-1+\dl}, \q \|(w^r,w^\th,w^z)\|_{L^\i(\Omega_\dl)}\leq C_\dl r^{-1+\dl},
\ee
where $\dl$ is suitably small and $\lim\limits_{\dl\rightarrow 0}C_\dl=+\i$.

\medskip
\subsection{Fast decay and vanishing of $u$}

\[
\]

{\bf Step 1.} almost first order decay of $|\na w |$.

Now we fix the number $\delta=0.1$.
From \eqref{e4.26} and \eqref{e4.27}, we can get
\be\label{4.41}
\|\big(\nabla w^r,\nabla w^z\big)\|_{L^2(\Omega_{\dl/2})}\leq C_\dl r^{-1+\dl},
\ee
and
\be\lab{4.42}
\|\nabla w^\th\|_{L^2(\Omega_{\dl/2})}\leq C_\dl r^{-1+\dl}.
\ee

Here it is more convenient to write the vorticity and  velocity in the Euclidean coordinates, which means
\be
w=w^1 (1,0, 0) +w^2(0, 1, 0) + w^3 (0, 0, 1) = w^r e_r+w^\th e_\th+w^ze_z \nn
\ee  for the vorticity and similar expression for the velocity.
From \eqref{4.41} and \eqref{4.42} we can obtain
\be\label{4.43}
\|\nabla w\|_{L^2(\Omega_{\dl/2})}\ls r^{-1+\dl}.
\ee
Now we fix any point $x=(r \cos \th, r \sin\th, z)$ where $r$ is large. Denote
\be
B_R=\{y \in \bR^3||x-y|\leq R\}.\nn
\ee

Next we introduce the following axially symmetric circling technique. We can find $r/2$ (round up to nearest integer) many balls, which are disjoint and generated by rotating $B_1$ around $z$ axis, such that their union is contained in the torus around the curve
\be
\big\{y|\s{y^2_1+y^2_2}=r,\ y_3=x_3\big\}.\nn
\ee
Since $u$ and $w$ is axially symmetric and the integrals of $|\na w|^2$ or $|\na u|^2$ on each ball is the same.

Then we can get
\be\label{4.44}
\begin{aligned}
&\int_{B_1}|\nabla w|^2dx\ls \f{1}{r}\int_{\Omega_{\dl/2}}|\nabla w|^2dx\ls r^{-3+2\dl},\\
&\int_{B_1}|\nabla u|^2dx\ls \f{1}{r}\int_{\Omega_{\dl/2}}|\nabla u|^2dx\ls r^{-1}.
\end{aligned}
\ee

In Euclid coordinates, the vorticity $w$ satisfy
\be\label{4.45}
-\Dl w=u\cdot \na w-(\na u)^T\cdot w,
\ee
where $\na u$ is the matrix
\be
\begin{bmatrix}
\p_1u^1 & \p_2u^1&\p_3 u^1\\
\p_1u^2 & \p_2u^2&\p_3 u^2\\
\p_1u^3 & \p_2u^3&\p_3 u^3
\end{bmatrix}\nn
\ee
and $(\na u)^T$ is the transpose of $\na u$.

Let $\p w$ be the derivative on $w$, which may represent $\p_1 w$, $\p_2 w$ or $\p_3w$. From \eqref{4.45}, we have
\be\label{4.46}
-\Dl \p w=\p u\cdot \na w+u\cdot \na \p w-(\na \p u)^T\cdot w-(\na u)^T\cdot \p w.
\ee

Now define a cut-off function $\phi(x)$ which is supported in $B_1$ and equal to $1$ in $B_{1/2}$. We can get
\be
-\Dl(\p w\phi)=\phi(-\Dl\p w)-2\na \p w\cdot \na \phi-\p w\Dl \phi.\nn
\ee
Then using the Green function $G(x^0,x)=\f{1}{4\pi}\f{1}{|x-x^0|}$ in three dimensions, we can get
\bea\lab{dw=i16}
\begin{aligned}
\p w(x)&=\int_{B_1}G\phi(-\Dl\p w) dx-2\int_{B_1}G\na \p w\cdot\na \phi dx-\int_{B_1}G\p w\Dl \phi dx \\
&=\underbrace{\int_{B_1}G\phi\p u\cdot \na wdx}_{I_1}+\underbrace{\int_{B_1}G\phi u\cdot \na \p wdx}_{I_2}\underbrace{-\int_{B_1}G\phi (\na \p u)^T\cdot wdx}_{I_3}\\
&\q\ \underbrace{-\int_{B_1}G\phi (\na u)^T\cdot \p wdx}_{I_4} \underbrace{-2\int_{B_1}G\na \p w\cdot\na \phi dx}_{I_5}\underbrace{-\int_{B_1}G\p w\Dl \phi dx}_{I_6}
\end{aligned}
\eea

Now we estimate $I_i(1\leq i\leq 6)$ term by term. From \eqref{4.44},
\be
|I_1|\ls \|\na w\|_{L^2(B_1)}\|\p u\|_{L^\infty(B_1)}\|G\|_{L^2(B_1)}\leq C_\dl r^{-3/2+\dl}.\nn
\ee
\be
|I_4|\ls \|\na u\|_{L^\i(B_1)}\|\p w\|_{L^2(B_1)}\|G\|_{L^2(B_1)}\leq C_\dl r^{-3/2+\dl}.\nn
\ee
\be
|I_6|\ls \|\Dl \phi\|_{L^\i(B_1)}\|\p w\|_{L^2(B_1)}\|G\|_{L^2(B_1)}\leq C_\dl r^{-3/2+\dl}.\nn
\ee
Using integration by parts, we have
\be\label{e4.47}
\begin{aligned}
|I_2|&\ls \Big|\int_{B_1}\p w\cdot(G\phi\na u+u\phi\na G+u G\na \phi)dx\Big|\\
&\leq C(\| u\|_{L^\i(B_1)}+\|\na u\|_{L^\i(B_1)})\|\p w\|_{L^2(B_1)}\|G\|_{L^2(B_1)}\\
&\q\ +\int_{B_1}\p w (u\phi\na G)dx\\
&\leq C_\dl r^{-3/2+\dl}+\|\p w\|_{L^\i(B_1)}\|u\|_{L^\i(B_1)}\int_{B_1}|\na G|\phi dx \\
&\leq C_\dl r^{-3/2+\dl}+C_\dl r^{-1+\dl}\\
&\leq C_\dl r^{-1+\dl}.
\end{aligned}
\ee Here we just used the usual boundedness  of $|\p w|$ and the decay of $u$ in \eqref{e4.42e}.

\be\label{e4.48}
\begin{aligned}
|I_3|&\ls \Big|\int_{B_1}(\p u)^T\cdot(w\phi\na G+G\na\phi w+ G\phi\na w)dx\Big|\\
 &\leq C\|\na w\|_{L^2(B_1)}\|G\|_{L^2(B_1)}\\
&\q\ +C(\|\p u\|_{L^\i(B_1)}\|w\|_{L^\i(B_1)})(\|\na G\|_{L^1(B_1)}+\| G\|_{L^1(B_1)})\\
&\leq C_\dl r^{-3/2+\dl}+\|\p u\|_{L^\i(B_1)}\|w\|_{L^\i(B_1)} \\
&\leq C_\dl r^{-3/2+\dl}+C_\dl r^{-1+\dl}\\
&\leq C_\dl r^{-1+\dl}.
\end{aligned}
\ee Here we just used the usual boundedness  of $|\p u|$ and the decay of $w$ in \eqref{e4.42e}.
Furthermore
\bea
|I_5|&\ls& \Big|\int_{B_1}(\p w)^T\cdot(\na G\cdot\na\phi+G\Dl \phi)dx\Big|\nn\\
&\leq&C\|\na w\|_{L^2(B_1)}(\|\na G\|_{L^2(B_1/B_{1/2})}+\| G\|_{L^2(B_1/B_{1/2})})\nn\\
&\leq& C_\dl r^{-3/2+\dl}.\nn
\eea

 At last, we can get, for large $r$
\be\label{4.47}
|\na w(r,z)| \leq C_\dl r^{-1+\dl}.
\ee

{\bf Step 2.} Almost $3/2$ order decay of $|\na w|$.
\medskip

First we need to obtain bounds on $|\na u|$.
Recall from \eqref{10.14} that
\[
\al
u^r(x)&=\int^{\pi}_{-\pi}\int^{2\pi}_0\int^\i_0 {G}(x,y)\psi \p_l w^\phi \cos \phi\rho d\rho d\phi dl\\
&\qq -2\int^{\pi}_{-\pi}\int^{2\pi}_0\int^\i_0 {G}(x,y)\p_\rho\psi\p_\rho u^\rho \cos\phi \rho d\rho d\phi dl\nn\\
&\qq -\int^{\pi}_{-\pi}\int^{2\pi}_0\int^\i_0 {G}(x,y)(\p^2_\rho\psi+\f{1}{\rho}\p_\rho\psi) u^\rho \cos\phi \rho d\rho d\phi dl.
\eal
\]Here we take $\psi$ as a standard cut-off function supported in the unit ball $B(x, 1)$ such that $\psi(x)=1$.
Differentiating this we obtain
\[
\al
\na u^r(x)&=\int^{\pi}_{-\pi}\int^{2\pi}_0\int^\i_0 \na {G}(x,y)\psi \p_l w^\phi \cos \phi\rho d\rho d\phi dl\\
&\qq -2\int^{\pi}_{-\pi}\int^{2\pi}_0\int^\i_0 \na {G}(x,y)\p_\rho\psi\p_\rho u^\rho \cos\phi \rho d\rho d\phi dl\nn\\
&\qq -\int^{\pi}_{-\pi}\int^{2\pi}_0\int^\i_0 \na {G}(x,y)(\p^2_\rho\psi+\f{1}{\rho}\p_\rho\psi) u^\rho \cos\phi \rho d\rho d\phi dl.
\eal
\]
 On  the first term of the righthand side, we use \eqref{4.47} and the bounds on $|\na G|$ in
Lemma \ref{l10.1}; on the second term, we use integration by parts to move the derivative on
 $u^\rho$ to other terms and then use the decay of $u^\rho$ in \eqref{e4.42e}.  Note this works since the singularity of the Green's function is cut off. The third term can be handled by the decay of $u^\rho$ and the bound on the gradient of the Green's function. Then we deduce that
\be\label{4.48}
|\p_r u^r(x)|+|\p_z u^r(x)|\leq C_\dl r^{-1+\dl}.
\ee

From this bound on $|\na u^r|$,  the divergence free condition
\[
\p_r u^r + \frac{u^r}{r} + \p_z u^z =0,
\]and the relation
\[
w^\th = \p_z u^r - \p_r u^z
\]we find that
\be
\lab{duz<}
|\p_r u^z|+|\p_z u^z|\leq C_\dl r^{-1+\dl}.
\ee Here we have used the decay of $w^\th$ in \eqref{e4.42e}.

Using the identities $w^r = -\p_z u^\th$ and $w^z = \p_r u^\th + \frac{1}{r} u^\th$, together
with \eqref{e4.42e}, we deduce
\be\label{4.53}
|\na u^\th|\leq C_\dl r^{-1+\dl}.
\ee
Combining \eqref{4.48},\eqref{duz<} and \eqref{4.53}, we arrive at
\be
\label{4.50}
|\na u|\leq C_\dl r^{-1+\dl}.
\ee Alternatively, one can also apply the mean value inequality in Theorem 1.7 \cite{Z:2} to derive this decay of $|\na u|$ from the decay of $u$.

Now we are ready to prove the almost $3/2$ order decay of $|\na w|$.
In  \eqref{dw=i16}, the formula for $\partial w$, only the term $I_2$ and $I_3$ need further analysis since all other terms decay at almost $3/2$ order already.

Inserting \eqref{4.47} into the third line of \eqref{e4.47}, we see that
\[
| I_2 | \le C_\dl r^{-3/2+\dl}.
\]Likewise, inserting \eqref{4.50} into the third line of \eqref{e4.48} yields
\[
| I_3 | \le C_\dl r^{-3/2+\dl}.
\]
Therefore,
\be
\lab{dw<-1.5}
|\p w|\leq  C_\dl r^{-3/2+\dl}.
\ee

{\bf Step 3.} Almost $3/2$ order decay of $|u|$ and vanishing.
\medskip

Now we come back to the formula of $u^r$ in \eqref{10.14}:
\[
\al
u^r(x)&=\int^{\pi}_{-\pi}\int^{2\pi}_0\int^\i_0 {G}(x,y)\psi \p_l w^\phi \cos \phi\rho d\rho d\phi dl\\
&\qq -2\int^{\pi}_{-\pi}\int^{2\pi}_0\int^\i_0 {G}(x,y)\p_\rho\psi\p_\rho u^\rho \cos\phi \rho d\rho d\phi dl\nn\\
&\qq -\int^{\pi}_{-\pi}\int^{2\pi}_0\int^\i_0 {G}(x,y)(\p^2_\rho\psi+\f{1}{\rho}\p_\rho\psi) u^\rho \cos\phi \rho d\rho d\phi dl.
\eal
\]This time we choose $\psi$ to be a standard cut-off function supported in $B(x, r/2)$ such that
$\psi=1$ in $B(x, r/4)$. From the estimate following \eqref{10.14}, the last two terms of the above identity already decays exponentially in $r$. Substituting \eqref{dw<-1.5} to the first term on the righthand side of the above identity and using the estimate of $G$, we see that
\be\label{4.56}
|u^r|+|\na u^r|\leq C_\dl r^{-3/2+\dl}.
\ee Here the gradient bound is obtained from differentiating the equation for $u^r$ and doing a similar estimate to that of $u^r$. Note that we are not claiming the decay of $|\na u^r|$ is better than that of $u^r$. So the proof is straight forward.

The bound \eqref{4.56} and
the divergence free condition
\[
\p_r u^r + \frac{u^r}{r} + \p_z u^z =0
\]imply that
\[
|\p_z u^z| \le C_\dl r^{-3/2+\dl}.
\] Moreover, \eqref{dw<-1.5} and the relation $w^r = - \p_z v^\th$ shows
\[
|\p_z u^\th| \le C_\dl r^{-3/2+\dl}.
\]
 By our assumption $\int^\pi_{-\pi}u^\th dz=\int^\pi_{-\pi}u^z dz=0$,  there exist, for any $r$, a $z_0, z_1$ such that $u^\th(r,z_0)=u^\th(r,z_1)=0$ where $z_0,z_1$ may be different for different $r$. Then
\be
|u^\th(r,z)|=|\int^z_{z_0}\p_z u^\th dz|\ls |\p_z u^\th(r, \cdot)|_{L^\i}2\pi\leq C_\dl r^{-3/2+\dl},
\ee
 and
\be
|u^z(r,z)|=|\int^z_{z_1}\p_z u^z dz|\ls |\p_z u^z(r, \cdot)|_{L^\i}2\pi\leq C_\dl r^{-3/2+\dl}.
\ee
So, all the above estimates imply that for large $r$
\be
|u(x)| \leq C_\dl  r^{-3/2+\dl}.
\ee
Recall that $\dl=0.1$;
 thus the decay rate of $u$ is faster than order $1$.  According to \cite{CSTY:1} or \cite{KNSS:1}, we have  proven that $u=0$, finishing the proof of Theorem \ref{thzp}. \qed

\medskip

Finally we prove Corollary \ref{codd}. Let $u$ be a D-solution given in the Corollary. Then $u$,
treated as functions of $r, z$,  satisfies:

\be\label{11.2}
\lt\{
\begin{aligned}
&(u^r\p_r+u^z\p_z)u^r-\f{(u^\th)^2}{r}+\p_r p=(\p^2_r+\f{1}{r}\p_r+\p^2_z-\f{1}{r^2})u^r,\\
&(u^r\p_r+u^z\p_z)u^\th+\f{u^ru^\th}{r}=(\p^2_r+\f{1}{r}\p_r+\p^2_z-\f{1}{r^2})u^\th,\\
&(u^r\p_r+u^z\p_z)u^z+\p_z p=(\p^2_r+\f{1}{r}\p_r+\p^2_z)u^z,\\
&\p_ru^r+\f{u^r}{r}+\p_zu^z=0,\\
&u^\theta |_{z=0, \pi}=0, \q  u^z |_{z=0, \pi}=0,\q  \p_z u^r |_{z=0, \pi}=0,
\end{aligned}
\rt.
\ee
where $x=(x',x_3)=(r\cos\th,r\sin\th,z)\in\bR^2\times[0,\pi]$. Now we extend the solution of \eqref{11.2} from $\bR^2\times [0,\pi]$ to $\bR^2\times[-\pi,\pi]$ by the following transformation. Define $\t{u}=\t{u}^r e_r+\t{u}^\th e_\th+\t{u}^z e_z$ in $\bR^2\times[-\pi,\pi]$ as follows
\be
\t{u}^r(r,z)=\lt\{
\begin{aligned}
&u^r(r,z)\qq\ \ \ (r,z)\in [0,\i]\times[0,\pi],\\
&u^r(r,-z)\qq (r,z)\in [0,\i]\times[-\pi,0],
\end{aligned}
\rt.\nn
\ee
\be
\t{u}^\th(r,z)=\lt\{
\begin{aligned}
&u^\th(r,z)\qq\q (r,z)\in [0,\i]\times[0,\pi],\\
&-u^\th(r,-z)\q (r,z)\in [0,\i]\times[-\pi,0],
\end{aligned}
\rt.\nn
\ee
\be
\t{u}^z(r,z)=\lt\{
\begin{aligned}
&u^z(r,z)\qq\q (r,z)\in [0,\i]\times[0,\pi],\\
&-u^z(r,-z)\q (r,z)\in [0,\i]\times[-\pi,0],
\end{aligned}
\rt.\nn
\ee
\be
\t{p}^z(r,z)=\lt\{
\begin{aligned}
&p^z(r,z)\qq\q (r,z)\in [0,\i]\times[0,\pi],\\
&p^z(r,-z)\q (r,z)\in [0,\i]\times[-\pi,0],
\end{aligned}
\rt.\nn
\ee
which mean that $u^r, p$ make even extension, while $u^\th,u^z$ make odd extension with respect to variable $z$. After such extension, it is easy to see that $(\t{u}^r,\t{u}^\th,\t{u}^z)$ satisfy \eqref{11.2} for $(r,z)$ in the domain $[0,+\i)\times[-\pi,\pi]$ with $p$ replaced by $\t{p}$. Since $\t{u}^\th$ and $\t{u}^z$ are odd with respect to $z$, we have
\be\label{11.4}
\int^\pi_{-\pi}\t{u}^\th(r,z)dz=\int^\pi_{-\pi}\t{u}^z(r, z)dz=0.
\ee

  Note that $\t{u}$ can be regarded as a $z-periodic$ solution on $\bR^2\times S^1$ with property \eqref{11.4}.  Theorem \ref{thzp} infers that $\t{u}=0$, hence $u=0$. \qed
\medskip

\remark By going through the proof of Theorem \ref{thzp}, one can see that the condition that
$u^z$, $u^\theta$ has $0$ mean in the $z$ direction can be relaxed to the following: for each $r>0$, there exists $z=z(r)$ such that $|u^z(r, z(r))|, |u^\theta(r, z(r))| \le \frac{C}{r}$.
Also, in case $u$ is in $L^2$ space, then the boundedness of the pressure $P$ will imply vanishing of $u$. We wish to thank Zhao Na for pointing out this fact to us.  In the 2 dimensional case, boundedness of $P$ is proven in \cite{GW:1}. In three dimensional exterior domains, Galdi  (\cite{Ga:1} Chapter X, Theorem 5.1), applies the decay properties of the Stokes kernel to prove rapid decays of the velocity, its gradient and second order derivative which implies the decay of $|\nabla P|$ from the equation. This in turn implies that $P$ converges to a constant at infinity. In general, using this method, it is not easy to prove that $P$ is bounded in the periodic case, due a lack of estimate of the Stokes kernel. Our decay estimates for $u$, $w$ and $\nabla w$ actually imply that $P$ can be chosen as a bounded function.

\medskip
{\bf Acknowledgement.}  We wish to thank Professors H. Brezis, H. J. Dong, Z. Lei and S.K. Weng
for helpful communications.  X. Pan is supported by Natural Science Foundation of Jiangsu Province (No.SBK2018041027) and National Natural Science Foundation of China (No. 11801268). Q. S. Zhang is grateful to the Simons' foundation for its support.
  Finally we are grateful to the two anonymous referees for the careful review and corrections and suggestions.

\bigskip

\noindent E-mail: \q B. Carrillo: \q bcarr011@ucr.edu; current address:
bcarrillo@saddleback.edu \\
\indent\qq\q \ X. Pan:\qq\q  xinghong\_87@nuaa.edu.cn\\
\indent\qq\q \    Q. S. Zhang: \ qizhang@math.ucr.edu


\begin{thebibliography}{00}

\bibitem[BG]{BG:1}Brezis, H.; Gallouet, T. {\it Nonlinear Schr\"odinger evolution equations.} Nonlinear Anal. 4 (1980), no. 4, 677-681.

\bibitem[Ch]{Ch:1}Chae, Dongho {\it Liouville-type theorems for the forced Euler equations and the Navier-Stokes equations.} Comm. Math. Phys.  326  (2014),  no. 1, 37-48.

\bibitem[CJ]{CJ:1} Choe, Hi Jun; Jin, Bum Ja {\it Asymptotic properties of axis-symmetric D-solutions of the Navier-Stokes equations.} J. Math. Fluid Mech. 11 (2009), no. 2, 208-232.

\bibitem[CPZZ]{CPZZ:1} Carrillo, B.; Pan, X.; Zhang, Q. S.; Zhao, N.  {\it Decay and vanishing of some {D}-solutions of the
  {N}avier-{S}tokes equations}. \newblock {arXiv:1808.10386}, 2018.

\bibitem[CW]{CW:1} Chae, Dongho; Wolf, J\"{o}rg {\it On Liouville type theorems for the steady Navier-Stokes equations in $\bR^3$.} J. Differential Equations  261  (2016),  no. 10, 5541-5560.

\bibitem[CSTY]{CSTY:1}
Chen, Chiun-Chuan;  Strain, Robert M.; Tsai, Tai-Peng; Yau, Horng-Tzer {\em Lower bound on the blow-up rate of the axisymmetric
 Navier-Stokes equations}. Int. Math Res. Notices (2008), vol. 8, artical ID rnn016, 31 pp.

\bibitem[Ga]{Ga:1}Galdi, G. P.  ``An introduction to the mathematical theory of the Navier-Stokes equations: Steady-state problems". Second edition. Springer Monographs in Mathematics. Springer, New York, 2011.

\bibitem[GW]{GW:1} Gilbarg, D.; Weinberger, H. F. {\it Asymptotic properties of steady plane solutions of the Navier-Stokes equations with bounded Dirichlet integral.} Ann. Scuola Norm. Sup. Pisa Cl. Sci. (4) 5 (1978), no. 2, 381-404.

\bibitem[Fr]{Fr:1961ACTA}  Finn, Robert {\it On the steady-state solutions of the Navier-Stokes equations, III}. Acta Math. 105 (1961), 197-244.

\bibitem[KNSS]{KNSS:1} Koch, G.; Nadirashvili  N.; Seregin, G.; Sverak, V.  {\em Liouville theorems for the Navier-Stokes equations and applications}. Acta Math. 203(2009), no. 1, 83-105.

\bibitem[KTW]{KTW:1} Kozono, Hideo; Terasawa, Yutaka; Wakasugi, Yuta {\it A remark on Liouville-type theorems for the stationary Navier-Stokes equations in three space dimensions.} J. Funct. Anal.  272  (2017),  no. 2, 804-818.

\bibitem[Le]{Le:1}   Leray, Jean, {\it \'Etude de diverses \'equations int\'egrales non lin\,eaires et de quelques probl\,emes que pose l'hydrodynamique.} (French) 1933. 82 pp.

\bibitem[LNZ]{LNZ:1} Lei, Zhen; Navas, Esteban A.; Zhang, Qi S. {\it A priori bound on the velocity in axially symmetric Navier-Stokes equations.} Comm. Math. Phys. 341 (2016), no. 1, 289-307.

\bibitem[LZ]{LZ:1} Lei, Zhen; Ren, Xiao; Zhang, Qi S.  {\it A Liouville Theorem for Axi-symmetric Navier-Stokes Equations on $R^2\times T^1$}. arXiv:1911.01571.

\bibitem[MB]{MB:2002CAMBRIDGE} Majda, Andrew J.; Bertozzi, Andrea L. ``Vorticity and incompressible flow". Cambridge Texts in Applied Mathematics, 27. Cambridge University Press, Cambridge, 2002. xii+545 pp. ISBN: 0-521-63057-6; 0-521-63948-4


\bibitem[Se]{Se:1} Seregin, G. {\it Liouville type theorem for stationary Navier-Stokes equations.} Nonlinearity  29  (2016),  no. 8, 2191-2195.

 \bibitem[TZ]{TZ:1}Tian, Gang; Zhang, Qi S. {\it Isoperimetric inequality under K\"{a}hler Ricci flow.} Amer. J. Math.  136  (2014),  no. 5, 1155-1173.

\bibitem[O]{O:1} O'Leary, Mike {\it Conditions for the local boundedness of solutions of the Navier-Stokes system in three dimensions}. Comm. Partial Differential Equations 28 (2003), no. 3-4, 617-636.

\bibitem[We]{We:1} Weng, Shangkun {\it Decay properties of smooth axially symmetric D-solutions to the steady Navier-Stokes equations}.  arXiv:1511.00752, J. Math. Fluid Mech. (2017).

\bibitem[Z1]{Z:1} Zhang, Qi S. {\it Global solutions of Navier-Stokes equations with large $L^2$ norms in a new function space.} Adv. Differential Equations 9 (2004), no. 5-6, 587-624.

\bibitem[Z2]{Z:2} Zhang, Qi S. {\it Local estimates on two linear parabolic equations with singular coefficients.} Pacific J. Math. 223 (2006), no. 2, 367-396.  Addendum arXiv:1905.13329.
\end{thebibliography}
\enddocument